%% file: main.tex
\documentclass[twoside]{article}

\usepackage{PRIMEarxiv}

\usepackage[utf8]{inputenc} 
\usepackage[T1]{fontenc}    
\usepackage{hyperref}       
\usepackage{url}            
\usepackage{booktabs}       
\usepackage{amsfonts}       
\usepackage{nicefrac}       
\usepackage{xcolor}         

\usepackage[numbers]{natbib}

\usepackage{amsmath}
\usepackage{amsfonts}
\usepackage{amssymb}
\usepackage{graphicx}
\usepackage{algorithm}
\usepackage{algpseudocode}
\usepackage{tikz}
\usepackage{xparse}
\usepackage{wrapfig}
\usepackage[bottom]{footmisc}
\usepackage{bbm}
\usepackage{multirow}
\usepackage{caption}
\usepackage{subcaption}
\usepackage{bm}
\usepackage{comment}
\usepackage{multibib}
\usepackage{enumitem}
\newcites{AP}{Additional References}

\newcommand{\va}{\mathbf{a}}
\newcommand{\vA}{\mathbf{A}}
\newcommand{\vb}{\mathbf{b}}
\newcommand{\vc}{\mathbf{c}}
\newcommand{\vC}{\mathbf{C}}

\newcommand{\vE}{\mathbf{E}}

\newcommand{\vI}{\mathbf{I}}
\newcommand{\vu}{\mathbf{u}}
\newcommand{\vU}{\mathbf{U}}
\newcommand{\vw}{\mathbf{w}}
\newcommand{\vW}{\mathbf{W}}
\newcommand{\vV}{\mathbf{V}}
\newcommand{\vx}{\mathbf{x}}
\newcommand{\vX}{\mathbf{X}}
\newcommand{\vh}{\mathbf{h}}

\newcommand{\vz}{\mathbf{z}}

\newcommand{\x}{x}

\newcommand{\cl}{\mathcal}

\newcommand{\vpi}{\boldsymbol{\pi}}
\newcommand{\valpha}{\boldsymbol{\alpha}}
\newcommand{\vbeta}{\boldsymbol{\beta}}

\newcommand{\vzero}{\boldsymbol{0}}
\newcommand{\vAdj}{\mathbf{A}^{\!(\text{adj})}}

\newcommand{\bR}{\mathbb{R}}
\newcommand{\bZ}{\mathbb{Z}}

\newcommand{\bone}{\mathbbm{1}}

\newcommand{\cI}{\mathcal{I}}
\newcommand{\cJ}{\mathcal{J}}

\newcommand{\cS}{\mathcal{S}}
\newcommand{\cV}{\mathcal{V}}
\newcommand{\cE}{\mathcal{E}}

\newcommand{\cN}{\mathcal{N}}
\newcommand{\cT}{\mathcal{T}}

\newcommand{\KL}{\text{KL}}
\newcommand{\Beta}{\text{Beta}}
\newcommand{\Ber}{\text{Ber}}

\newcommand{\Dz}{D_{\text{z}}}
\newcommand{\Dx}{D_{\text{x}}}
\newcommand{\Dc}{D_{\text{c}}}

\newcommand{\vzb}{\mathbf{z}^{(\text{b})}}
\newcommand{\vzc}{\mathbf{z}^{(\text{c})}}
\newcommand{\vAb}{\mathbf{A}^{(\text{b})}}
\newcommand{\vAc}{\mathbf{A}^{(\text{c})}}
\newcommand{\vcc}{\mathbf{c}^{(\text{c})}}
\newcommand{\vcb}{\mathbf{c}^{(\text{b})}}

\newcommand{\prs}{\texttt{routing-sm}}
\newcommand{\pr}{\texttt{routing}}
\newcommand{\fl}{\texttt{facility-loc}}
\newcommand{\tsp}{\texttt{tsp}}

\newcommand{\rems}{\texttt{revenue-max-sm}}
\newcommand{\rem}{\texttt{revenue-max}}
\newcommand{\eg}{\texttt{energy-grid}}
\newcommand{\ca}{\texttt{caching}}

\newcommand{\ours}{\texttt{MIPnet}}
\newcommand{\ourssup}{\texttt{MIPnet-sup}}
\newcommand{\oursuns}{\texttt{MIPnet-unsup}}
\newcommand{\nd}{\texttt{N-div}}
\newcommand{\todo}[1]{}
\renewcommand{\todo}[1]{{\color{red} TODO: {#1}}}

\title{Solving Recurrent MIPs with Semi-supervised Graph Neural Networks
}

\author{
  Konstantinos Benidis \\
  Amazon Research \\
  Berlin, Germany\\
  \texttt{kbenidis@amazon.com} \\
   \And
   Ugo Rosolia \\
   Amazon Research \\
  Luxembourg\\
  \texttt{urosolia@amazon.com} \\
  \And
  Syama Rangapuram \\
  Amazon Research \\
  Berlin, Germany\\
  \texttt{rangapur@amazon.com} \\
  \AND
George Iosifidis\thanks{Work done while at Amazon.}\\
  Delft University of Technology  \\
  Delft, Netherlands\\
  \texttt{g.iosifidis@tudelft.nl} \\
  \And
Georgios Paschos\\
  Amazon Research \\
  Luxembourg\\
  \texttt{gpaschos@amazon.com} \\
}

\date{}

\begin{document}
\maketitle

\begin{abstract}
We propose an ML-based model that automates and expedites the solution of MIPs by predicting the values of variables. Our approach is motivated by the observation that many problem instances share salient features and solution structures since they differ only in few (time-varying) parameters. Examples include transportation and routing problems where decisions need to be re-optimized whenever commodity volumes or link costs change. Our method is the first to exploit the sequential nature of the instances being solved periodically, and can be trained with ``unlabeled'' instances, when exact solutions are unavailable, in a semi-supervised setting. Also, we provide a principled way of transforming the probabilistic predictions into integral solutions. Using a battery of experiments with representative binary MIPs, we show the gains of our model over other ML-based optimization approaches.

\end{abstract}


\input{texfiles/intro.tex}
\input{texfiles/related.tex}

\input{texfiles/model_v2.tex}

\input{texfiles/experiments.tex}

\input{texfiles/conclusion.tex}

\bibliography{mlopt_updated.bib}

\input{texfiles/supplement.tex}

\bibliographystyleAP{plainnat}
\bibliographyAP{mlopt_bib_app.bib}

\end{document}

%% file: texfiles/intro.tex

\section{Introduction}
\label{sec:intro}
	
The solution of mixed integer programming (MIP) problems is as important as it is challenging to achieve. Indeed, MIP is employed for optimizing real-world operations, decision processes and systems, including transportation  \cite{laporte09, johnson05}, facility location \cite{farahani09}, production planning \cite{pochet2006, chen10}, content distribution networks \cite{paschos-fnt}, all the way to structured predictions \cite{vivek2012, kulesza08} and neural network (NN) training \cite{leng-aaai18}. However, despite the success and wide availability of MIP solvers, the fast and accurate optimization of such problems remains largely elusive.
This has given rise to perpetual efforts for approximation algorithms \cite{shmoys-book} and heuristic solutions \cite{boussaid-heuristics}, and for improving the performance of solvers \cite{bixby}. 

Nonetheless, an aspect that has received less attention is that several MIP instances share properties and common structures. In fact, practitioners, more often than not, need to solve repeatedly problems that differ only in few parameters, cf. \cite{misra2019, bertsimas2021voice}. For example, routing and other transportation problems are re-solved periodically over the same graph whenever new demands or costs are obtained. This motivates the use of ML in order to explore correlations among problem properties and solution values, which in turn can be leveraged to expedite the solution of new instances, enabling their solution even in real-time as an increasing number of applications require \cite{bertsimas2019online}. The idea of ML-assisted optimization is not new per se; it has been successfully applied to configure solvers \cite{kleinberg-neurips19}, and design heuristics, e.g., in Branch-and-Bound (BnB) techniques \cite{gasse2019exact}. Importantly, recent studies focused on learning problem structures towards \emph{predicting} (some) variables  \cite{gasse2019exact, dingaaai20, nair2020solving, karalias202Nips} and/or active constraints \cite{bertsimas2019online, bertsimas2021voice}, with promising results -- we review them in Sec. \ref{sec:related}.

In this work we make several steps towards the latter direction by developing an ML-based method for facilitating the solution of \emph{binary} MIPs which (might) exhibit temporal structure. For example, internet traffic follows a diurnal (and not random) pattern \cite{diotsigm04, torino-dataset}; and the same holds for a plenitude of transportation problems \cite{hidalgo-2008}. Motivated by these observations, we propose \ours{}, a Long Short-Term Memory (LSTM)-based probabilistic variable-prediction model \cite{lstm97} that operates on a permutation and scale-invariant embedding space created by a Graph Convolutional Network (GCN) \cite{battaglia2018relational}. The GCN exploits the intuitive bipartite graph MIP representation \cite{gasse2019exact, nair2020solving}, towards encoding the problems' salient properties without manual feature engineering. 

Furthermore, unlike prior works, \ours{} employs \emph{semi-supervised} learning to augment the training data without computational costs, while still benefiting from the robustness of supervised learning when possible. And including this unsupervised loss provides additional gains. Namely, not all variables have the same impact on the MIP's objective and on constraint violation. Indeed, misprediction of some variables can render the problem infeasible, while others affect only the objective (and to different extent). Clearly, while a model using only supervised loss cannot discern such conditions, \ours{} can learn these effects with per-variable granularity, and mitigate their impact accordingly. Our goal is to predict a significant portion of the binary variables so as to facilitate the optimization of the remaining (binary and continuous) variables. The selection of these variables is based on a confidence metric following a principled Bayesian approach that enables tunable variable selection without (expensive) model re-training.

\ours{} is evaluated through a battery of experiments, using both real and synthetic traces, and in diverse real-world problems that include network routing, facility location and TSP, among others. These are representative problems in operational research, and have been used in prior works \cite{bertsimas2019online, bertsimas2021voice, nair2020solving, dingaaai20} and references therein. The experiments reveal significant gains over these previous works in terms of variable prediction accuracy, while maintaining high constraint feasibility. We also demonstrate the importance of exploiting the temporal dimension of these problems (whenever it is prevalent), towards improving our predictions. Finally, we stress that we focus on binary MILPs which have applications in a vast range of real-world problems;\footnote{E.g., 164 of 240 benchmarks in \cite{miplib} are binary/integer LPs; while $90\%$ of variables are binary in the rest.} yet, our method paves the road towards tackling general MIPs.

In summary, the contributions of this work are:
\vspace{-2mm}
\begin{itemize}[leftmargin=2mm]
    \item We propose \ours{}, a GCN-based probabilistic model that learns to solve MIPs and selects the most confident variables based on a principled Bayesian approach. \ours{} speeds up significantly the solution of \emph{general binary} MIPs while exploiting, for the first time, any temporal relation across the problem instances.

\vspace{-1mm}

    \item \ours{} employs semi-supervised learning to benefit from labeled data when available, and remain operative when these are hard to obtain. This hybrid approach further allows the model to identify and prioritize variables that are pivotal for feasibility and optimality.

\vspace{-1mm}    
    
    \item In a series of experiments, \ours{} is proved to find accurate solutions to a range of different problems, consistently outperforming SotA benchmarks such as \cite{bertsimas2019online, nair2020solving}.
\end{itemize}

\textbf{Paper Structure}. Sec. \ref{sec:related} reviews the related work and Sec. \ref{sec:model} introduces the  model, training and inference approach. Sec. \ref{sec:experiments} presents the experimental results and Sec. \ref{sec:conclusions} concludes the study. Model and dataset details, and additional experiments are provided as supplementary material.

%% file: texfiles/related.tex

\vspace{-1mm}
\section{Related work}\label{sec:related}
\vspace{-1mm}

\paragraph{Learning Configurations \& Heuristics} A first thrust of works learn to tune (hyper) parameters of solvers \cite{hutter09, hutter11, lopez16, kleinberg-neurips19, ansotegui-aistats15, balcan-branch-18}. For example, \cite{hutter09, hutter11} leverage local search to identify good configurations and \cite{lopez16} developed a pertinent software library. These suggestions are oblivious to, and hence unable to benefit from, the problem structure. 
A different approach is to learn heuristics. Many works in this context learn how to select BnB variables \cite{land-bnb, khalil-branch-16,gasse2019exact, alvarez-branching-17, hansknecht-branch18} or BnB nodes \cite{he-neurips14, song-arxiv18}, and how to run primal heuristics \cite{khalil-branch-17}. Similar ideas are used in cutting planes \cite{balcan2021sample} and in optimization over graphs \cite{dai-icml16, dai2017learning, kim2021learning}. On the other hand, \cite{bonami18, kruber17} (and references therein) learn optimization models via alternative decompositions and problem reformulations.

\paragraph{Identifying Substructures} The third thrust of related work focuses on learning and exploiting problem substructures. Early efforts include predicting \emph{backdoor variables} \cite{dilkina09-16, dilkina09-17}, i.e., instantiating key variables to increase the problem's tractability. Similarly, assuming the availability of training data, \cite{misra2019, klauco2019} propose sampling methods for predicting active constraints. Along these lines, \cite{bertsimas2019online, bertsimas2021voice} builds Optimal Classification Trees and Feedforward NNs to predict active constraints and variables. Similarly, \cite{xavier-informs21} assigns variables using a $k$-Nearest Neighbours classifier. Finally, \cite{song2020general, wu2021learning} use Imitation Learning and RL for splitting ILPs. These works do not employ richer (and more promising) problem representations, nor they design bespoke NN architectures for the problems at hand.

\paragraph{Employing GCNs} To that end, \cite{gasse2019exact, nair2020solving, dingaaai20} use GCNs where the MIPs are encoded with bipartite graphs. This approach is permutation-invariant, thus can yield {generalizable} learned policies. In \cite{gasse2019exact} the GCN is used for variable selection in BnB; \cite{nair2020solving} learns also how to initialize (\emph{Neural Diving}); and \cite{dingaaai20} uses cross-entropy loss to predict probability distributions for the binary variables so as to expedite branch exploration. On the other hand, \cite{karalias202Nips} proposes a GNN-based unsupervised approach for graph-related problems, where the integral solution is recovered with derandomization through sequential decoding. This interesting approach yields feasible solutions with high probability, yet is not applicable to general MIPs.

Unlike these works (see overview \cite{bengio2021machine}), \ours{} considers the \emph{temporal evolution} of problem's parameters; employs semi-supervised learning to reduce the requirement for training datasets; and follows a Bayesian approach to variable selection through learned distributions with tunable confidence and infeasibility penalization. Appendix \ref{app:related} includes further comparisons with prior work.

%% file: texfiles/model_v2.tex


\vspace{-1mm}
\section{Semi-supervised Temporal GCN}\label{sec:model}
\vspace{-1mm}

\subsection{Problem Setting}

Consider the binary MILP in standard form:\footnote{We do not specifically include equality constraints since they can be rewritten as two inequality constraints.} 
\begin{equation}
	\begin{aligned} 
		\underset{\vz}{\text{minimize}} & \quad \vc^T \vz\\
		\text{subject to} & \quad \vA \vz  \preceq \vb, \\
	\end{aligned}
	\label{eq:binary_MILP} 
\end{equation}
where $\vz \!= [\vzb ; \vzc]$, with $\vzb \!\in \{0, 1\}^{\Dz^{(\text{b})}}$, and $\vzc \!\in \bR^{\Dz^{(\text{c})}}$ being the binary and continuous variables, respectively, and $\Dz \!=\! \Dz^{(\text{b})} \! + \Dz^{(\text{c})}$. An instance of \eqref{eq:binary_MILP} is defined by the set of parameters $\phi \!= \! \{\vc, \vA, \vb \}$, with $\vc \! \in \bR^{\Dz}$, $\vA \!\in \bR^{\Dc \times\Dz}$, $\vb \! \in\bR^{\Dc}$ (notation details in Appendix \ref{app:notation}).

We focus on applications where \eqref{eq:binary_MILP} needs to be solved recurrently with different parameters that follow a temporal structure. Consider a set of time series $\cS$ with $\cT_s$ timesteps, where each timestep corresponds to an instance with parameters $\phi_{s, t}$, with $s\in \cS$ and  $t\in\cT_s$.  
For the  set of  instances $\{\phi_{s,t}\}_{s\in \cS, t\in \cT_s}$, consider their optimal solutions $\{\vz^\star_{s, t}\}_{s\in \bar{\cS}, t\in \bar{\cT}_s}$, with $\bar{\cS}\subseteq \cS$ and $\bar{\cT}_s \subseteq \cT_s$, i.e., we assume the availability of solution only for a subset of instances.
Our goal is to train a model that learns a globally shared mapping from the problem parameters $\phi_{s, t}$ to a probabilistic representation $p(\vz^\star_{s, t})$ of the optimal solution.
Leveraging this mapping, we present a method to select a subset of variables for which the optimal assignment is known with high confidence -- the size of the subset is a user-defined parameter. Then, we fix these variables in problem~\eqref{eq:binary_MILP} and we solve the resulting lower dimensional sub-problem.
As shown in the result section, 
for every new set of instances $\{\phi_{s,t}'\}_{s\in \cS', t\in \cT_s'}$, 
the proposed methods is able to expedite the computation of the optimal solution.

We elaborate next on the model architecture, training and inference procedures. For clarity of exposition we drop the time series and timestep subscripts unless necessary.

\subsection{Proposed Method}

Inspired by \cite{gasse2019exact, dingaaai20, nair2020solving}, we represent an instance of \eqref{eq:binary_MILP} as a bipartite graph by having one node for each variable and constraint, and connecting a variable node to a constraint node iff the variable is active on that constraint. We use a GCN block to map the nodes to an embedding space, followed by an LSTM block to capture the temporal structure of the MILP instances $\{\phi_{s,t}\}_{s\in \cS, t\in \cT_s}$. The model is trained in a semi-supervised manner, exploiting both a supervised and an unsupervised component.

Two key properties of \eqref{eq:binary_MILP} are the scale and permutation invariance, i.e., the solution of the problem does not change if we scale (appropriately) the problem, or permute its parameters and variables. If a model does not satisfy these properties, it would need to be re-trained for each variation of a problem. Our architecture satisfies both scale invariance by introducing a normalization step and permutation invariance by allowing only symmetric functions among nodes.

\paragraph{Parameter Normalization}

We initially normalize all parameters $\phi$ of \eqref{eq:binary_MILP} in order to make the training more stable and the model scale-invariant. This means that the model will learn to solve the problem with parameters drawn from a ``base'' distribution, and every scaled problem version can be solved just by normalizing the parameters. The normalization needs to ensure that the relative weights of the objective parameters and each constraint remain the same; therefore we normalize them separately. The normalized parameters $\forall i \leq \Dc, j \leq \Dz$, are computed as follows:
\begin{equation}\label{eq:param_norm}
a_{i,j}\! =\!  \frac{a_{i, j}}{\| [\va_i^T ; b_i] \|_p},  \ b_i \! =  \! \frac{b_i}{\| [\va_i^T ; b_i]\|_p}, \ c_j \!= \!  \frac{c_j}{\|\vc\|_p},
\end{equation}
where $\va_i$ is the $i$-th row of $\vA$. We use $p=2$. For problem instances of different sizes the normalization can lead to a distribution shift. We discuss this in Appendix \ref{app:distr_shift}.

\paragraph{Model Architecture.}
The model builds on the bipartite graph representation of MILPs. Formally, consider a graph $G = (\cV, \cE, \vAdj)$ defined by the set of nodes $\cV$, with $|\cV| = \Dz + \Dc$, the set of edges $\cE$, with $|\vE|$ equal to the number of nonzero entries in the constraint matrix $\vA$, and the graph adjacency matrix $\vAdj = [\vI_{\Dz}, \vA^T ; \vA,  \vI_{\Dc}]\in \bR^{|\cV|  \times |\cV|}$, 
i.e., the (non-binary) adjacency matrix contains the coefficients from the constraint matrix $\vA$ and self-loops for all nodes given by the identity matrices. One set of $\Dz$ nodes in the bipartite graph corresponds to variables and the other set of $ \Dc$ nodes to constraints. An edge $e_{i,j}\in\cE$ indicates that variable $j$ appears in the $i$-th constraint.

The nodes of the graph are associated with a set of features $\vU = \vU(\phi) \in\bR^{|\cV| \times D_u}$ derived by the problem parameters $\phi$, where $D_u$ is the feature dimension. The feature vector of each node is constructed by linearly combining in a permutation invariant way all the variable and constraint parameters that the node is related to.
We refer the reader to Appendix \ref{app:feature} for a detailed description.

The graph representation of \eqref{eq:binary_MILP} motivates using a GCN \cite{kipf2016semi}. A GCN with $L$ layers is denoted: 
\begin{equation}
\label{eq:GNN_embedding}
\vX^{(l+1)} \!=\! g(\vX^{(l)}, \vAdj ; \theta_g^{(l)}), \ \forall \ l = 0, \ldots, L-1,
\end{equation}
with $\vX^{(0)} = \vU$ and $\theta_g^{(l)}\in \Theta$ the learnable parameters of each layer. The adjacency matrix defines the graph connectivity and determines how information is aggregated at each layer, while the number of layers define the number of hops away that a node gets information from.  Note that the GCN is applied to each node in parallel and the resulting embedding is by construction permutation invariant.
We follow the GCN propagation rule as defined in \cite{kipf2016semi, battaglia2018relational}, with a few (optional) modifications between each layer of the GCN as in \cite{nair2020solving}: (i) We include a nonlinearity after each layer, (ii) we include skip connection inputs and (iii) we apply layer normalization \cite{ba2016layer} in the output of each layer.

The output of the final layer $\vX = \vX^{(L)} \in \bR^{|\cV|\times \Dx}$ is the final embedding of each node which will be used as input to the next NN block of the model architecture. We use an LSTM network to capture the temporal evolution of \eqref{eq:binary_MILP}. The inputs to the LSTM at time $t$ are the embedded nodes  $\vX$ as well as the previous network output $\vh_{s, t-1}$, i.e., $
\vh_{s,t} = h(\vh_{s,t-1}, \vX_{s,t-1} ; \theta_h)$,
with $\theta_h\in \Theta$ the learnable parameters of the LSTM. The LSTM retains the permutation invariance since it acts on the feature dimension. The network output at time $t$, i.e., $\vh_{s,t} $, is appropriately projected to the parameters of the selected probabilistic representation of each variable using a multilayer perceptron (MLP). In general the projection has the following form:
\begin{equation}
\label{eq:model_out}
\psi_{s, t, j} = f(\vx_{s, t, j}, \vh_{s, t} ; \theta_f), 
\end{equation}
where $\theta_f\!\in\! \Theta$ are the learnable MLP parameters and $\vx_{s, t, j} \!=\! [\vX_{s, t}]_j$ are the embedded features of the $j$-th variable node that acts as a skip connection, while the constraint nodes are discarded. Note that the mapping layer can apply any function in the feature dimension of $\vh_{s, t}$ but only symmetric functions across the node dimension. Here, $\psi$ represents a generic set of parameters for the binary variables that may differ based on the selected model output. If \eqref{eq:binary_MILP} includes continuous variables the model outputs their values and not a distribution, i.e.,  $\psi_{s, t}\! =\! \{p(\vzb_{s,t}), \hat{\vz}^{\text{(c)}}_{s,t} \}$. In particular, we leverage two MLPs to project the output of the LSTM to the estimates of binary and continuous variables (for a discussion about the prediction of continuous variables please see Appendix \ref{app:continuous}).
The probabilistic representation of binary variables becomes concrete in Sec.  \ref{ssec:supervised}.

\subsection{Training}

A common characteristic of MILP applications is the difficulty in obtaining labels. Specifically, to obtain one label, one has to solve an instance of the MILP, which depending on the use case may take from minutes to hours or even days. It is therefore of high interest for the generality of \ours~to be able to train with a semi-supervised loss, exploiting labeled data and benefiting from unlabeled data.


\begin{figure*}[t]
\centering
	\subfloat{\scalebox{0.65}{\input{figures/model.tikz}}}
	\caption{\ours{} --  At each time step $t$, we construct the features $\vU_{s,t}$ of all nodes based on the MILP temporal parameters, feed them into a GCN that produces the final embedding $\vX_{s,t}$. The embeddings are fed to an  LSTM whose output is  projected to the set of parameters $\psi_{s,t}$ and the corresponding loss $\ell_{s,t}$ is computed.}
	\label{fig:temporal_supervised}
\end{figure*}

\paragraph{Supervised Setting.}
\label{ssec:supervised}

For the instances that have available labels (optimal solutions) we can train the model with an MLE approach. A natural choice for a probabilistic representation of binary variables is the Bernoulli distribution $\Ber(z ; \pi)$ \cite{nair2020solving, wu2021learning}, where $\pi$ is the Bernoulli parameter. Given a Bernoulli representation of a binary variable it is straightforward to assign a value to the variable, but hard to quantify how accurate or confident this prediction is (see Appendix \ref{app:bernoulli} for a discussion). A different approach was proposed in \cite{nair2020solving} that used the SelectiveNet \cite{geifman2019selectivenet} and learned two sets of binary variables. One set indicates if a variable is going to be selected and the other its binary value (to be used if selected). Here, we introduce a principled way to define how accurate or confident is a prediction via its own variance. We leverage this notion of variance and follow a Bayesian approach by using the Beta distribution $\Beta(\pi ; \alpha, \beta)$ to model the Bernoulli parameters, where $\alpha$ and $\beta$ are the parameters of the Beta distribution. This approach can be readily used as a statistically principled selection method, and is a key advantage of our method. 

Now, we elaborate on the model output \eqref{eq:model_out}. The model yields the Beta distribution parameters $\psi \!=\! \{\alpha, \beta\}$ (we ignore the continuous variables and drop subscripts), and for each variable with optimal value $z^\star$ the likelihood is computed as $
\int_{0}^{1} \Beta(\pi ; \alpha, \beta) \Ber(z^\star ; \pi) d\pi$, i.e., we need to integrate over all values of $\pi$. Unfortunately this integral does not have a closed-form solution; yet, we can approximate it efficiently using a quadrature method (Monte-Carlo sampling can be also used, but is not as efficient). We use the Clenshaw-Curtis quadrature method \cite{clenshaw1960method} where the function to be integrated is evaluated at the $K$ roots of a Chebyshev polynomial and the integral is approximated by a weighted average of the integral's function $I(\pi)\! =\! \Beta(\pi ; \alpha, \beta) \Ber(z^\star ; \pi)$ in specific and predefined points  $\{\bar{\pi}_k\}_{k=0}^{K/2}$, i.e., 
\begin{equation}
\int_{0}^{1} \!\Beta(\pi ; \alpha, \beta) \Ber(z^\star ; \pi) d\pi \!=\! \int_{0}^{1} \! I(\pi) d\pi \approx \vw^T I(\bar{\vpi}), \notag
\end{equation}
where $\vw\in\mathbb{R}^{K/2+1}$ is independent of the function (precomputed) and $I(\bar{\vpi}) = [I(\bar{\pi}_0), \dots, I(\bar{\pi}_{K/2})]^T$. For details of the method we refer the reader to Appendix \ref{app:ccq}.
With the above formulation the negative log-likelihood (NLL) of all the instances becomes:
\begin{equation}
\label{eq:suploss}
\ell_{\text{sup}}(\vz^\star, \psi)
\approx - \sum_{s=1}^{|\cS|}\sum_{t=1}^{|\cT_s|} \sum_{j=1}^{\Dz} \log\left(\vw^T I(\bar{\vpi}_{s, t, j}) \right),   
\end{equation}
where $\psi\!=\!\{\valpha_{s,t}, \vbeta_{s,t} \}$.

In Appendix \ref{app:reg} we propose an additional regularization term, while in Appendix \ref{app:weight} we describe a weighted version of \ref{eq:suploss} for imbalanced labels.

\paragraph{Unsupervised Setting.}

Apart from  \eqref{eq:suploss}  we consider an unsupervised loss that is optimized together and is essentially the objective of \eqref{eq:binary_MILP} with a penalty term for the constraint violations:
\begin{equation}
\label{eq:unsup_loss}
\ell_{\text{unsup}} = \sum_{s=1}^{|\cS|}\sum_{t=1}^{|\cT_s|} \vc_{s,t}^T \hat{\vz}_{s,t} + \lambda_{\text{c}} \|\left(\vA_{s,t}  \hat{\vz}_{s,t}  - \vb_{s,t} \right)_+ \|^2,
\end{equation}
where $\lambda_{\text{c}} \geq 0$ and $(x)_+ = \max(x, 0)$. Here, $\hat{\vz}_{s,t}$ represents a prediction based on the learned Beta distribution that is computed as 
$\hat{z} = \sigma\left(\frac{\alpha}{\alpha + \beta}\right)$,
where $\sigma(\cdot)$ is the sigmoid function and $\frac{\alpha}{\alpha + \beta}$ is the mean of Beta. Applying the sigmoid to the mean gives values close to the limits that we use as a proxy for the true binary since exact rounding would block the gradient flow. Further, in the unsupervised loss we include the continuous variables that are directly given by the network output \eqref{eq:model_out}, i.e., $\psi_{s, t} =$ $\{\valpha_{s, t}, \vbeta_{s, t}, \hat{\vz}^{(\text{c})}_{s, t}\}$ where $\valpha_{s, t}, \vbeta_{s, t} \in \bR^{\Dz^{(\text{b})}}$ and
$\hat{\vz}^{(\text{c})}_{s, t} \in \bR^{\Dz^{(\text{c})}}$.


The unsupervised loss, as explained earlier, allow us to work with instances for which labels are hard to obtain and, at the same time, provides the means to learn the impact of erroneous predictions (with per-variable granularity) on the objective function and constraint violation. This, in turn, is proved an effective practical control for mitigating the effect of the (few) miss-predictions that \ours{} yields. Note also that the continuous variables are optimized jointly with the binary through this loss. 

Putting the above together, the overall loss is given by
\begin{equation}
\label{eq:loss}
\ell = \ell_{\text{sup}} + \lambda \ell_{\text{unsup}},
\end{equation}
with $\lambda \!\geq\! 0$. Note that if a label is not available, the supervised loss of that instance is masked and we compute only the unsupervised loss. Fig. \ref{fig:temporal_supervised} summarizes the architecture.

\subsection{Inference and Variable Selection}
Our prediction problem is very challenging since it can return an infeasible solution even if it mispredicts only few variables. Hence, instead of targeting an all-or-nothing solution (as in \cite{bertsimas2021voice}, see Table \ref{tab:bertsimas}), we apply a \emph{variable selection} method where we generate a recommended assignment for each variable associated with a notion of confidence. By fixing the value of high-confidence variables and calling an MILP solver on the reduced problem, we speed up the discovery of high-quality solutions, as shown also in \cite{dingaaai20, nair2020solving, bertsimas2021voice, xavier-informs21}. Differently to these works however, we take a formal Bayesian approach into this variable selection problem, which additionally, does not require retraining.

Given a trained model and a new instance $\phi' = \{\vc, \vA, \vb \}$ the model outputs the parameters $\psi' = \{\valpha, \vbeta, \hat{\vz}^{(\text{c})}\}$ with $\valpha, \vbeta \in \bR^{\Dz^{(\text{b})}}$ and $\hat{\vz}^{(\text{c})} \in \bR^{\Dz^{(\text{c})}}$. For each of the $\Dz^{(\text{b})}$ binary variables we compute the mean and the standard deviation of its Beta distribution, given by $\mu_j = \frac{\alpha_j}{\alpha_j  \!+ \beta_j}$ and $\sigma_j = \sqrt{\frac{\alpha_j \beta_j}{(\alpha_j + \beta_j)(\alpha_j + \beta_j +1)}}$.
Since a confident output has mean close to 0 or 1 and low variance, we define the score
$s_j = \min(\mu_j, 1 - \mu_j) + \gamma \sigma_j$,  for $j = 1, \dots, \Dz^{(\text{b})}$,
and we select a desired percentage $\rho$ of the variables with the lowest score. Parameter $\gamma \!\geq\! 0$ is tuned on a validation set. For a given set of selected variables, we fix their values rounding their mean, and optimize the remaining binary and continuous variables.
Following this approach we solve the remaining problem only once, although one has the option to take multiple samples from the learned distribution, solve the remaining problems multiple times and keep the best solution \cite{nair2020solving}.

%% file: figures/model.tikz
\begin{tikzpicture}[shorten >=1pt,->,draw=black!60, font=\footnotesize, scale=0.5]%
			\definecolor{lilac}{HTML}{e4cbff}
			\tikzstyle{every edge}=[draw,-stealth, thin];
			\tikzstyle{box}=[rectangle, draw==black!50, fill=lilac!40, minimum height=12pt, inner sep=0pt]
			\tikzstyle{input node}=[circle, draw, thick, minimum size=30pt, inner sep=2pt];
			\tikzstyle{embed node}=[box, minimum width=30pt, inner sep=2pt];
			\tikzstyle{gcn node}=[box, minimum width=40pt, inner sep=2pt];
			\tikzstyle{lstm node}=[box, minimum width=30pt, inner sep=2pt];
			\tikzstyle{map node}=[box, minimum width=30pt, inner sep=2pt];
			\tikzstyle{loss node}=[circle, draw, thick, minimum size=30pt, inner sep=2pt];
			
			\tikzstyle{control node}=[circle, draw, thick, fill=lilac!40, minimum width=25pt, inner sep=2pt];
			\tikzstyle{count node}=[circle, draw, thick, minimum size=25pt, inner sep=2pt];
			\tikzstyle{switch node}=[circle, draw, thick, minimum size=25pt, inner sep=2pt];
			\tikzstyle{state node}=[circle, draw, thick, minimum size=25pt, inner sep=2pt];
			\tikzstyle{observation node}=[circle, draw, thick, fill=lilac!40, minimum width=25pt, inner sep=2pt];
			
			\foreach \i/\j in {1/{t-1}, 2/{t}, 3/{t+1}} {%
				\node[loss node] (l\i) at (10*\i,11.9) {$\ell_{s,\j}$};%
				\node[map node] (psi\i) at (10*\i,9.4) {$\psi_{s,\j} = \{p(\vzb_{s,\j}), \hat{\vz}^{\text{(c)}}_{s,\j} \}$};%
				\node[lstm node] (h\i) at (10*\i,7.1) {$\vh_{s,\j}$};%
				\node[gcn node] (x\i) at (10*\i,4.8) {$\vX_{s,\j} ; \vAdj_{s,\j}$};%
				\node[embed node] (u\i) at (10*\i,2.5) {$\vU_{s, \j}$};%
				\node[input node] (phi\i) at (10*\i,0) {$\phi_{s, \j}$};%
				\path (phi\i) edge (u\i);%
				\path (u\i) edge (x\i);%
				\path (x\i) edge (h\i);%
				\path (h\i) edge (psi\i);%
				\path (psi\i) edge (l\i);%
			}%
			
			\node [left of=l1, node distance=1.2in, align=flush right] (ltlabel) {Loss};%
		        \node [left of=psi1, node distance=1.2in, align=flush right] (psitlabel) {Projection};%
   			\node [left of=h1, node distance=1.2in, align=flush right] (htlabel) {LSTM};%
		        \node [left of=x1, node distance=1.2in, align=flush right] (xtlabel) {GCN};%
 			\node [left of=u1, node distance=1.2in, align=flush right] (utlabel) {Features};%
   			\node [left of=phi1, node distance=1.2in, align=flush right] (phitlabel) {Input};%
			
			\foreach \i in {1, 2, 3} {%
				\path (x\i) edge [out=150, in=215] (psi\i);%
			}%
			
			\foreach \i/\j in {1/2, 2/3} {%
				\path (h\i) edge (h\j) ;%
			}%
			
\end{tikzpicture}

%% file: texfiles/experiments.tex

\section{Experiments} \label{sec:experiments}
\vspace{-2mm}
We evaluate \ours{} using carefully-selected problems and datasets, as detailed in Appendix \ref{app:datasets}. For the $\pr$ problem we used the dataset {\tt{Geant}} \cite{SNDlib10} collecting real-world \emph{temporal} traffic data from Internet Service Providers; for the \fl{} we used the real-world topology from the {\tt{nobel-germany}} dataset \cite{SNDlib10}; and for the \ca{} problem we used actual requests from the {\tt{MovieLens}} datasets \cite{movielens}. We generated synthetic data for the \tsp{}, \eg{} and \rem{} problems. We ran every experiment 3 times in order to evaluate model variance. Appendices \ref{app:exp_details} and \ref{app:add_exp} include evaluation details and additional experiments, respectively.

\subsection{Evaluating Accuracy \& Feasibility} \label{sec:exp_acc}

We assess the accuracy, optimality and feasibility of \ours{} against various state of the art methods. Namely, we compare our results with Neural Diving (\nd) \cite{nair2020solving}, a GCN model that selects variables using the SelectiveNet approach \cite{geifman2019selectivenet}; the method proposed in \cite{bertsimas2021voice}, an approach that learns a classifier from problem parameters to complete solutions (based on the unique solutions appeared in the training set); and a variant of our model trained only with supervised loss, \ourssup{}, to evaluate the importance of the unsupervised loss component. As a further ablation, we also evaluate our method when trained only with unsupervised loss which, however, did not achieve satisfactory performance (see Appendix \ref{app:add_exp}). As a final method we use the MIP solver SCIP \cite{scip}. For all problems expect for \pr, SCIP finds an optimal solution for all instances within the $15$-minutes time limit. On the other hand for \pr, SCIP finds an optimal solution for $99.73\%$ of instances within the $60$-minutes time limit, which highlights the problem's complexity. These optimal solutions -- and the suboptimal ones when optimals are not available  -- serve as \textit{ground truth}, and the results of all other methods are presented relative to them.

\begin{table*}[t]
	\footnotesize
	\center
	\caption{\small Accuracy (mean \textpm\ std)$\%$ of \ours{} vs. \nd{}. Bold indicates the best method (higher values are better).}
	\resizebox{0.68\textwidth}{!}{\begin{tabular}{cc | cccccc}
		\toprule
		\midrule
		   $\rho \ (\%)$  & Method   &       \pr &          \fl &           \tsp  &        \eg &        \rem   &  \ca\\
		\midrule
		\midrule
		\multirow{3}{*}{30\%} & \ours &  \textbf{97.60 \textpm\ 0.13} &  \textbf{100.00 \textpm\ 0.00} &  \textbf{98.12 \textpm\ 1.38} &  \textbf{100.00 \textpm\ 0.00} &  \textbf{100.00 \textpm\ 0.00} & \textbf{100.00 \textpm\ 0.00}\\
		& \ourssup &  \textbf{97.60 \textpm\ 0.04} &  99.88 \textpm\ 0.16 & 94.56 \textpm\ 1.73 &  99.97 \textpm\ 0.05  &  \textbf{100.00 \textpm\ 0.00} & \textbf{100.00 \textpm\ 0.00}\\
		& \nd &  84.40 \textpm\ 0.72 &  \textbf{100.00 \textpm\ 0.00 }& 93.43 \textpm\ 0.00 &  \textbf{100.00\textpm\ 0.00}  &  \textbf{100.00 \textpm\ 0.00} & \textbf{100.00 \textpm\ 0.00}\\
		
		\midrule
		
		\multirow{3}{*}{40\%} & \ours &  \textbf{96.03 \textpm\ 0.13} &  \textbf{100.00 \textpm\ 0.00} &  \textbf{97.54 \textpm\ 1.73} &  \textbf{100.00 \textpm\ 0.00} &  \textbf{100.00 \textpm\ 0.00} & \textbf{100.00 \textpm\ 0.00}\\
		& \ourssup &  96.01 \textpm\ 0.15 &  99.87 \textpm\ 0.16  & 93.97 \textpm\ 1.49 &  99.97 \textpm\ 0.04  &  \textbf{100.00 \textpm\ 0.00} & \textbf{100.00 \textpm\ 0.00}\\
		& \nd &  87.17 \textpm\ 3.68 &  99.98 \textpm\ 0.02 &  96.09 \textpm\ 1.76 &  99.20 \textpm\ 1.13  &  \textbf{100.00 \textpm\ 0.00} & \textbf{100.00 \textpm\ 0.00}\\
		
		\midrule

		\multirow{3}{*}{50\%} & \ours &  \textbf{94.46 \textpm\ 0.09} &  \textbf{100.00 \textpm\ 0.00} &  \textbf{96.97 \textpm\ 1.94} &  \textbf{99.99\textpm\ 0.00} &  \textbf{100.00 \textpm\ 0.00} & \textbf{100.00 \textpm\ 0.00}\\
		& \ourssup &  \textbf{94.46 \textpm\ 0.11} &  99.79 \textpm\ 0.24 & 93.42 \textpm\ 1.40 &  99.94 \textpm\ 0.03  &  \textbf{100.00 \textpm\ 0.00} & \textbf{100.00 \textpm\ 0.00}\\
		& \nd &  89.34 \textpm\ 0.11 &  99.91 \textpm\ 0.03 &  95.02 \textpm\ 0.54 &  99.61 \textpm\ 0.53 &  \textbf{100.00 \textpm\ 0.00} & \textbf{100.00 \textpm\ 0.00}\\
		
		\midrule
		
		\multirow{3}{*}{60\%} & \ours &  92.83 \textpm\ 0.08 &  \textbf{99.98 \textpm\ 0.01} &  \textbf{96.38 \textpm\ 2.05} &  \textbf{99.96 \textpm\ 0.01} &  \textbf{100.00 \textpm\ 0.00} & \textbf{100.00 \textpm\ 0.00}\\
		& \ourssup &  \textbf{93.15 \textpm\ 0.11} &  99.71 \textpm\ 0.29 & 92.96 \textpm\ 1.35 &  99.86 \textpm\ 0.13  &  \textbf{100.00 \textpm\ 0.00} & \textbf{100.00 \textpm\ 0.00}\\
		& \nd &  90.38 \textpm\ 2.55 &  99.91 \textpm\ 0.02 &  95.41 \textpm\ 2.61 &  \textbf{99.96 \textpm\ 0.01} &  \textbf{100.00 \textpm\ 0.00} & \textbf{100.00 \textpm\ 0.00}\\
	
		\midrule

		\multirow{3}{*}{70\%} & \ours &  91.25 \textpm\ 0.07 &  \textbf{99.84 \textpm\ 0.05} &  \textbf{95.72 \textpm\ 2.05} &  \textbf{99.80 \textpm\ 0.01} &  \textbf{100.00 \textpm\ 0.00} & \textbf{100.00 \textpm\ 0.00}\\
		& \ourssup &  \textbf{91.86 \textpm\ 0.09} &  99.41 \textpm\ 0.58 & 92.62 \textpm\ 1.19 &  99.74 \textpm\ 0.14   &  \textbf{100.00 \textpm\ 0.00} & \textbf{100.00 \textpm\ 0.00}\\
		& \nd &  90.30 \textpm\ 0.92 &  99.64 \textpm\ 0.06 &  93.93 \textpm\ 0.84 &  99.64 \textpm\ 0.11 &  \textbf{100.00 \textpm\ 0.00} & \textbf{100.00 \textpm\ 0.00}\\
	
		\bottomrule
	\end{tabular}}
	\label{tab:accuracy}
\end{table*}

\begin{table*}[t]
	\footnotesize
	\center
	\caption{\small Infeasibility (mean \textpm\ std)$\%$ of \ours{} vs. \nd{}. Bold indicates the best method (lower values are better).}
	\resizebox{0.68\textwidth}{!}{\begin{tabular}{cc | cccccc}
		\toprule
		\midrule
		 $\rho \ (\%)$    &  Method  &       \pr &          \fl &           \tsp  &        \eg &        \rem  &  \ca \\
		\midrule
		\midrule
		\multirow{3}{*}{30\%} & \ours &   \textbf{59.09 \textpm\ 3.31} &   \textbf{0.00 \textpm\ 0.00} &  \textbf{0.09 \textpm\ 0.13} &  \textbf{0.00 \textpm\ 0.00} &  \textbf{0.00 \textpm\ 0.00} & \textbf{0.00 \textpm\ 0.00}\\
		& \ourssup &  61.27 \textpm\ 0.04 &   \textbf{0.00 \textpm\ 0.00} &  75.52 \textpm\ 38.37 &  0.02 \textpm\ 0.03 &  \textbf{0.00 \textpm\ 0.00} & \textbf{0.00 \textpm\ 0.00}\\
		& \nd &  100.00 \textpm\ 0.00 &   \textbf{0.00 \textpm\ 0.00} &  31.22 \textpm\ 17.52 &  \textbf{0.00 \textpm\ 0.00} &  \textbf{0.00 \textpm\ 0.00} & \textbf{0.00 \textpm\ 0.00}\\
		
		\midrule
		
		\multirow{3}{*}{40\%} & \ours &   \textbf{99.97 \textpm\ 0.04} &   \textbf{0.02 \textpm\ 0.01} &  0.46 \textpm\ 0.58 &  \textbf{0.00 \textpm\ 0.00} &  \textbf{0.00 \textpm\ 0.00} & \textbf{0.00 \textpm\ 0.00}\\
		& \ourssup &  \textbf{99.97 \textpm\ 0.04} &   0.11 \textpm\ 0.04 &  80.79 \textpm\ 27.16 &  0.02 \textpm\ 0.03 &  \textbf{0.00 \textpm\ 0.00} & \textbf{0.00 \textpm\ 0.00}\\
		& \nd &  100.00 \textpm\ 0.00 &  7.82 \textpm\ 11.06 &  \textbf{0.18\textpm\ 0.26} & \textbf{0.00 \textpm\ 0.00} &  \textbf{0.00 \textpm\ 0.00} & \textbf{0.00 \textpm\ 0.00}\\
		
		\midrule

		\multirow{3}{*}{50\%} & \ours &  \textbf{100.00 \textpm\ 0.00} &  \textbf{0.07 \textpm\ 0.05} &  \textbf{2.83 \textpm\ 2.78} &  \textbf{0.00 \textpm\ 0.00} & \textbf{0.00 \textpm\ 0.00} & \textbf{0.00 \textpm\ 0.00}\\
		& \ourssup &  \textbf{100.00 \textpm\ 0.00} &   1.03 \textpm\ 0.03 &  94.26 \textpm\ 8.12 &  0.09 \textpm\ 0.08 &  \textbf{0.00 \textpm\ 0.00} & \textbf{0.00 \textpm\ 0.00}\\
		& \nd & \textbf{100.00 \textpm\ 0.00} &  0.28 \textpm\ 0.40 & 76.14 \textpm\ 8.49 &  0.09 \textpm\ 0.08 &  \textbf{0.00 \textpm\ 0.00} & \textbf{0.00 \textpm\ 0.00}\\
		
		\midrule
		
		\multirow{3}{*}{60\%} & \ours &  \textbf{100.00 \textpm\ 0.00} &  \textbf{0.24 \textpm\ 0.18} &  \textbf{13.06 \textpm\ 7.82} &  \textbf{0.04 \textpm\ 0.03} & \textbf{0.00 \textpm\ 0.00} & \textbf{0.00 \textpm\ 0.00} \\
		& \ourssup &  \textbf{100.00 \textpm\ 0.00} &   3.78 \textpm\ 0.75 &  99.8 \textpm\ 0.28 &  0.09 \textpm\ 0.08 &  \textbf{0.00 \textpm\ 0.00} & \textbf{0.00 \textpm\ 0.00}\\
		& \nd & \textbf{100.00 \textpm\ 0.00} &  5.19 \textpm\ 4.19 &  34.07 \textpm\ 46.62 &  \textbf{0.04 \textpm\ 0.03} &  \textbf{0.00 \textpm\ 0.00} & \textbf{0.00 \textpm\ 0.00} \\

		\midrule

		\multirow{3}{*}{70\%} & \ours &  \textbf{100.00 \textpm\ 0.00} &  \textbf{2.92 \textpm\ 3.47} &  48.90 \textpm\ 17.19 &  \textbf{0.02 \textpm\ 0.03} &  \textbf{0.00 \textpm\ 0.00} & \textbf{0.00 \textpm\ 0.00}\\
		& \ourssup &  \textbf{100.00 \textpm\ 0.00} &   10.83 \textpm\ 0.39 &  100.00 \textpm\ 0.00 &  0.11 \textpm\ 0.11 &  \textbf{0.00 \textpm\ 0.00} & \textbf{0.00 \textpm\ 0.00}\\
		& \nd &  \textbf{100.00 \textpm\ 0.00} &  19.75 \textpm\ 12.90 &  \textbf{48.23 \textpm\ 36.64} &  0.04 \textpm\ 0.06 &  \textbf{0.00 \textpm\ 0.00} & \textbf{0.00 \textpm\ 0.00}\\

		\bottomrule
	\end{tabular}}
	\label{tab:infeasibility}
\end{table*}

\begin{table*}[t]
	\footnotesize
	\center
	\caption{\small  Optimality gap (mean \textpm\ std)$\%$ of \ours{} vs. \nd{}. Bold indicates the best method (lower values are better).}
	\resizebox{0.68\textwidth}{!}{\begin{tabular}{cc | cccccc}
			\toprule
			\midrule
			$\rho \ (\%)$    &  Method  &       \pr &          \fl &           \tsp  &        \eg &        \rem  &  \ca\\
			\midrule
			\midrule
			\multirow{3}{*}{30\%} & \ours &  6.59 \textpm\ 0.15 &  0.07 \textpm\ 0.01 &  \textbf{0.24 \textpm\ 0.02} &  \textbf{0.00 \textpm\ 0.00} &  \textbf{0.00 \textpm\ 0.00} & \textbf{0.00 \textpm\ 0.00}\\
			& \ourssup\ &  \textbf{6.03 \textpm\ 0.63} &  \textbf{0.04 \textpm\ 0.05} &  2.38 \textpm\ 0.00 &  \textbf{0.00 \textpm\ 0.00} &  \textbf{0.00 \textpm\ 0.00} & \textbf{0.00 \textpm\ 0.00}\\
			& \nd\ &  - &  0.28 \textpm\ 0.23 &  0.25 \textpm\ 0.04 &  \textbf{0.00 \textpm\ 0.00} &  \textbf{0.00 \textpm\ 0.00} & \textbf{0.00 \textpm\ 0.00}\\
			
			\midrule
			
			\multirow{3}{*}{40\%} & \ours &  \textbf{6.27 \textpm\ 0.00} &  \textbf{0.16 \textpm\ 0.01} &  0.37 \textpm\ 0.07 &  \textbf{0.00 \textpm\ 0.00} &  \textbf{0.00 \textpm\ 0.00} & \textbf{0.00 \textpm\ 0.00}\\
			& \ourssup\ &  27.59 \textpm\ 0.00 &  0.20 \textpm\ 0.10 &  4.20 \textpm\ 0.00 &  \textbf{0.00 \textpm\ 0.00} &  \textbf{0.00 \textpm\ 0.00} & \textbf{0.00 \textpm\ 0.00}\\
			& \nd\ & -  &  0.44 \textpm\ 0.35 &  \textbf{0.26 \textpm\ 0.04} &  0.08 \textpm\ 0.04 &  \textbf{0.00 \textpm\ 0.00} & \textbf{0.00 \textpm\ 0.00}\\
			
			\midrule
			
			\multirow{3}{*}{50\%} & \ours &  - &  \textbf{0.34 \textpm\ 0.03} &  \textbf{0.74 \textpm\ 0.19} &  \textbf{0.00 \textpm\ 0.00} &  \textbf{0.00 \textpm\ 0.00} & \textbf{0.00 \textpm\ 0.00}\\
			& \ourssup\ &  - &  0.72 \textpm\ 0.42 &  6.05 \textpm\ 0.00 &  \textbf{0.00 \textpm\ 0.00} &  \textbf{0.00 \textpm\ 0.00} & \textbf{0.00 \textpm\ 0.00}\\
			& \nd\ &  -  &  3.55 \textpm\ 4.88 &  5.70 \textpm\ 0.63 &  8.81 \textpm\ 12.45 &  \textbf{0.00 \textpm\ 0.00} & \textbf{0.00 \textpm\ 0.00} \\
			
			\midrule
			
			\multirow{3}{*}{60\%} & \ours &  - &  \textbf{0.73 \textpm\ 0.10} &  1.49 \textpm\ 0.30 &  \textbf{0.00 \textpm\ 0.00} &  \textbf{0.00 \textpm\ 0.00} & \textbf{0.00 \textpm\ 0.00}\\
			& \ourssup\ &  - &  2.15 \textpm\ 1.72 &  4.18 \textpm\ 0.00 &  0.01 \textpm\ 0.01 &  \textbf{0.00 \textpm\ 0.00} & \textbf{0.00 \textpm\ 0.00}\\
			& \nd\ &  - &  1.22 \textpm\ 0.55 &  \textbf{0.48 \textpm\ 0.01} &  0.01 \textpm\ 0.00 &  \textbf{0.00 \textpm\ 0.00} & \textbf{0.00 \textpm\ 0.00}\\
			
			\midrule
			
			\multirow{3}{*}{70\%} & \ours &  -  &  \textbf{2.66 \textpm\ 0.72} &  2.08 \textpm\ 0.17 &  \textbf{0.01 \textpm\ 0.00} &  \textbf{0.00 \textpm\ 0.00} & \textbf{0.00 \textpm\ 0.00}\\
			& \ourssup\ &  - &  5.80 \textpm\ 5.04 &  - &  \textbf{0.01 \textpm\ 0.01} &  \textbf{0.00 \textpm\ 0.00} & \textbf{0.00 \textpm\ 0.00}\\
			& \nd\ &  - &  5.83 \textpm\ 1.97 & \textbf{1.79 \textpm\ 0.21} &  0.04 \textpm\ 0.03 &  \textbf{0.00 \textpm\ 0.00} &  \textbf{0.00 \textpm\ 0.00}\\
			
			\bottomrule
	\end{tabular}}
	\label{tab:opt_gap}
\end{table*}

For a problem with $\Dz^{(b)}$ binary variables, we fix the $\big \lceil{\rho\Dz^{(b)}}\big \rceil$ most-confident ones, where $\rho\in(0,1)$, and use a solver for the exact optimization of the remaining subproblem, subject to \ours{}'s assignments. Table~\ref{tab:accuracy} summarizes the percentage of correct predictions on the assigned variables (\emph{accuracy}) across all instances, while Table~\ref{tab:infeasibility} presents the percentage of infeasible instances. Table~\ref{tab:opt_gap} compares the optimality gap of the different approaches, which is computed using \textit{only} the set of feasible solutions (denoted with $\cN_f$), and is defined as the average difference between a method's objective (\texttt{obj}) and SCIP's optimal value (\texttt{obj}$^*$), i.e., $\texttt{opt\_gap} \!=\! \frac{1}{|\cN_f|} \sum_{\texttt{obj} \in \cN_f} (\texttt{obj}\! -\! \texttt{obj}^*)$. We repeat the experiments for different $\rho$ values.\footnote{In \ours{}, we can decide a different portion of assigned variables without model retraining. This stands in stark contrast to \nd{} where the variable selection percentage is predefined.} 

We observe that \ours{} matches or outperforms \nd{} in almost all cases, while it exhibits more consistent behaviour across runs. Compared to \ourssup{}, it is either on par or better, especially in feasibility.
A notable observation is that accuracy itself is not adequate to validate the performance of a model. For example, in the \pr\ and \tsp\ datasets we observe that accuracy values well-above $90\%$ can lead to (almost) $100\%$ infeasibility; while in \fl\ accuracy values $>99\%$ can yield a significant portion of infeasible instances -- especially for \ourssup{} and \nd{}. This is associated with the complexity of the constraints and highlights the importance of using an unsupervised component in model training. 

\begin{table}[t]
	\center
	\caption{\small Accuracy, infeasibility and optimality gap of \cite{bertsimas2021voice}.}
	\scriptsize
	{\begin{tabular}{c | c |c | c}
			\toprule
			 & \textbf{Accuracy} & \textbf{Infeasibility} & \textbf{Optimality gap} \\
			\midrule
            \pr & - & - & - \\
            \fl & 98.13 & 0.00 & 1.78  \\
            \tsp & 86.29 & 0.00  & 19.61 \\
            \rem & 99.22 & 100.00 & - \\
            \ca & 99.91  & 0.00 & 52.49 \\
			\bottomrule
		\end{tabular}
	}
	\label{tab:bertsimas}
\end{table}

\begin{figure*}[t]
	\centering
	\begin{subfigure}[b]{0.30\textwidth}
		\centering
		\includegraphics[width=0.8\textwidth]{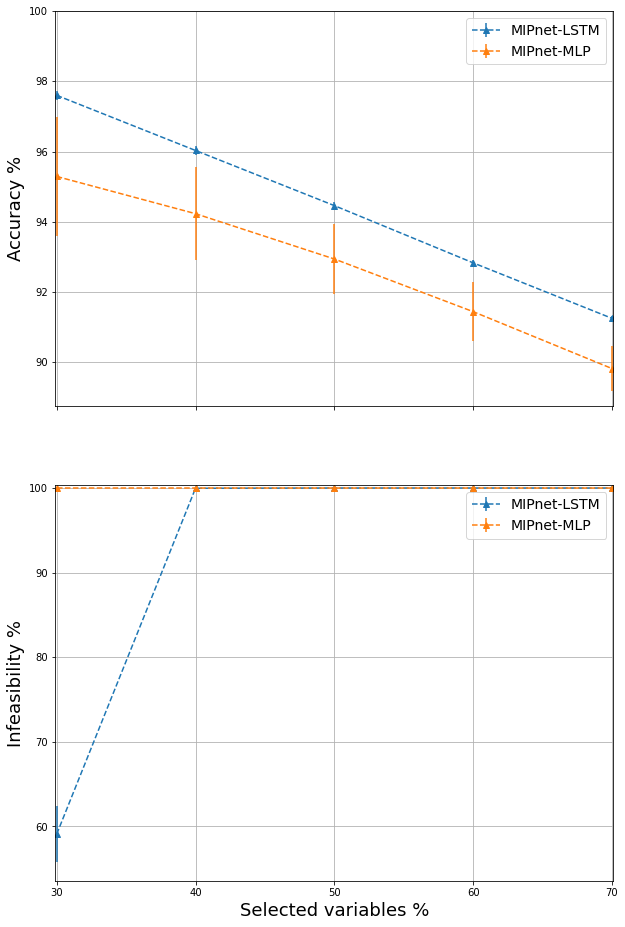}
		\caption{\pr}
		\label{fig:pr_acc}
	\end{subfigure}
	\begin{subfigure}[b]{0.30\textwidth}
		\centering
		\includegraphics[width=0.82\textwidth]{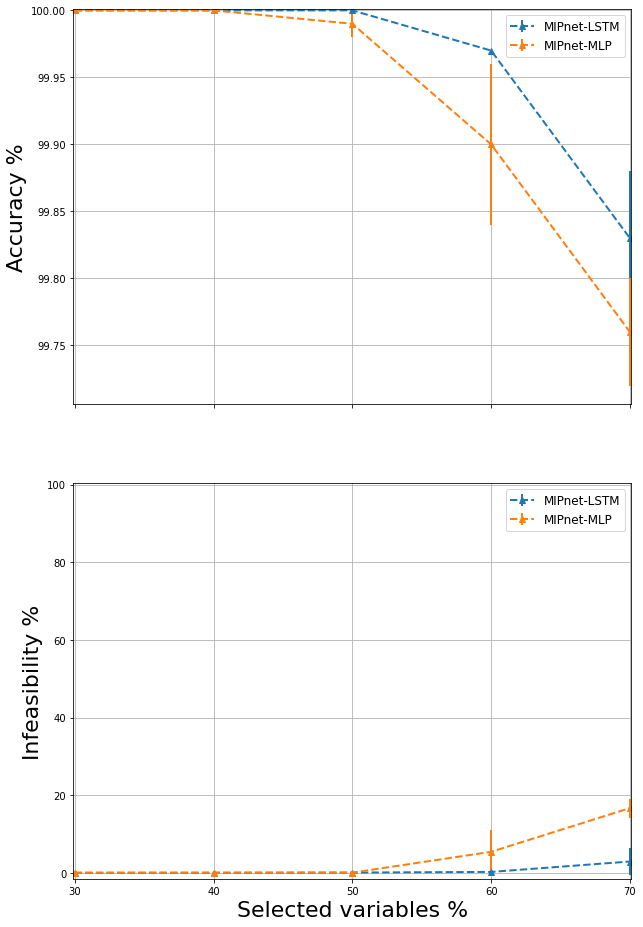}
		\caption{\fl}
		\label{fig:fl_acc}
	\end{subfigure}
	\begin{subfigure}[b]{0.30\textwidth}
		\centering
		\includegraphics[width=0.8\textwidth]{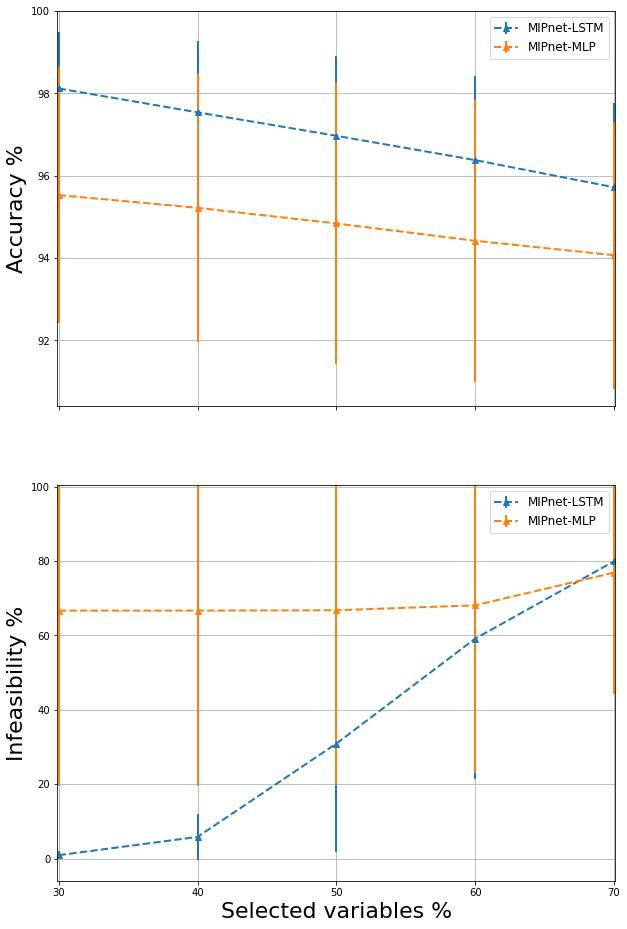}
		\caption{\tsp}
		\label{fig:tsp_acc}
	\end{subfigure}	\caption{Comparison of \ours\ w\textbackslash\ and w\textbackslash o a temporal component.}
	\label{fig:temporal}
\end{figure*}

Since \cite{bertsimas2021voice} is an approach that selects a complete solution, i.e., there is no option for selecting $\rho\%$ of variables, we present the results of this method separately. 
As shown in Table~\ref{tab:bertsimas}, this all-or-nothing approach learns feasible solutions for \tsp, \fl, and \ca\ (no variance was observed across runs). In these problems the constraints are time-invariant (do not change across instances). Hence, optimal solutions from the training dataset are feasible for the test dataset, and indeed we see that ~\cite{bertsimas2021voice} finds feasible solutions. Alas, these are often far from optimal, e.g., for \tsp{} and \ca{} the average gap is $19.61 \%$ and $52.49\%$, respectively. On the other hand, for \rem{} where the constraints change across instances, optimal solutions from the training dataset are infeasible for new problem instances and the method returns infeasible variable assignments. That is, this approach does not generalize to solutions not seen in the training dataset. Finally, we highlight that we were able to run the approach from~\cite{bertsimas2021voice} only on a smaller version of the \pr{} dataset, which alludes to potential implementation difficulties of this approach -- see Appendix \ref{app:add_exp}. 


\subsection{Evaluating the Temporal Component of \ours{}}In this section we explore the benefit of leveraging the temporal aspect of the problem instances. Figure \ref{fig:temporal} illustrates the differences in accuracy and infeasibility when \ours{} employs an MLP instead of an LSTM, for three problems with different temporal properties. Namely, the real-world demands in \pr{} follow a diurnal pattern; the demands in \fl{} have a strong time-dimension; while the edge costs in \tsp{} are generated by adding random perturbations in successive instances.

We observe that the LSTM-based \ours{} outperforms consistently, and most-often substantially, its MLP-based version both in terms of mean and variance values. In other words, all else being equal in the model and inference approach, the addition of LSTM has substantial gains, without requiring any compromises (MLP never outperforms LSTM). This finding underlines the importance of leveraging the temporal evolution pattern of the problem parameters, a hitherto overlooked aspect in all prior works. 

\begin{figure}
	\centering
	\includegraphics[width=0.95\columnwidth]{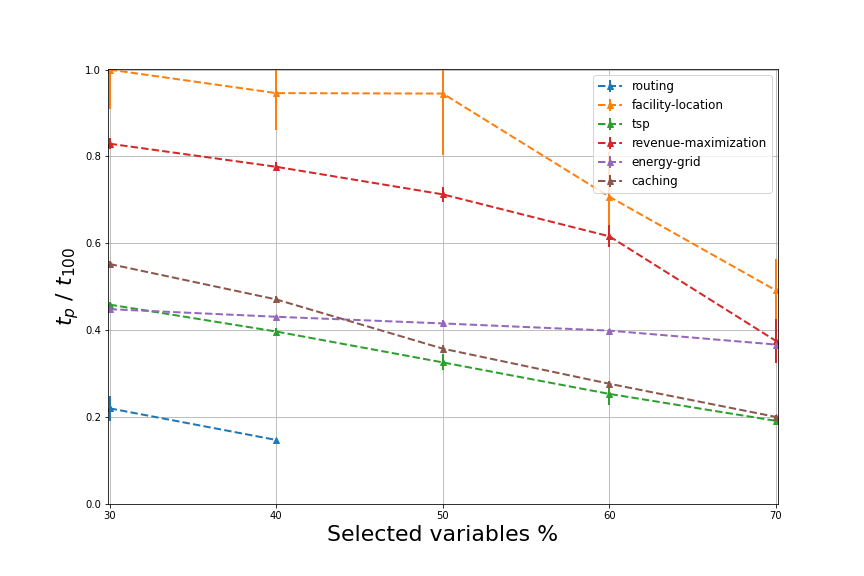}
	\vspace{-3mm}
	\caption{Calibrated running time for various $\rho$ values. $t_p$ is the solution time when a percentage $\rho$ of variables is assigned by \ours{} and $p\!=\!1\!-\!\rho$ variables are assigned by the solver. $t_{100}$ is the time when using only the solver.}
	\label{fig:time}
\end{figure}

\subsection{Trading Off Accuracy and Solution Speed via $\rho$}
\vspace{-2mm}
Finally, we explore how the percentage of assigned variables $\rho$, affects the solution time. These experiments reveal also the gains of \ours{} in terms of solution speed, compared to off-the-shelf solvers. In detail, we denote with $t_p$ the time the SCIP solver requires for solving the optimization problem, where $p = 1 - \rho$ is the percentage of variables that are \emph{not} assigned by \ours{} (hence, need to be optimized by the solver). Figure~\ref{fig:time} reports the normalized time defined as the ration $t_p / t_{100}$, i.e., the ratio between the solver (SCIP) time when $\rho$ variables are assigned using \ours{}, and the time for solving the entire problem with SCIP.\footnote{These times are calibrated in order to account for the different computational environments we have used in the experiments.} As expected, the time decreases when a higher percentage of variables is assigned. For all problems an average of, at least, $2\times$ speedup is obtained by setting $\rho = 0.7$, while for \rem, \eg, and \ca{} we achieve a 2x speedup already for $\rho = 0.3$, i.e., when assigning only $30\%$ of the variables, and a maximum $5 \times$ speed up. These experiments manifest the solution speed gains of \ours{} compared to a solver, and the importance of our variable selection threshold which indeed can be used to trade off accuracy with solution speed. 

%% file: texfiles/conclusion.tex
\section{Conclusions}
\label{sec:conclusions}

Solving MIPs via ML methods can revolutionize the solution of large-scale NP-hard problems and avail MIP to new application domains. This potential is hampered by unprecedented challenges: the values of variables are inherently correlated and need to be jointly predicted; even tiny value assignment errors may lead to infeasibility or unbounded suboptimality; training data are not readily available; and generalizing trained models is highly non-trivial. In this work, we propose \ours{}, a systematic framework for tackling general binary MILPs that may exhibit temporal structure, which: employs GCNs to automate embeddings and capture correlation among variables; uses an LSTM to account for, hitherto-overlooked, time dependencies across problem instances (a regularly-encountered operational condition); and includes a principled Bayesian assignment strategy with user-tunable confidence. We follow a semi-supervised approach so as to account for the lack of solved instances, a practical limitation when dealing with NP-hard problems, and in order, also, to identify key variables impacting feasibility and optimality. Our work fills (some of) the gaps in the literature and, we believe, contributes towards establishing ML as an off-the-shelf MIP technology, including for non-binary and nonlinear problems that were not addressed by the proposed approach. 

%% file: texfiles/supplement.tex
\onecolumn

\centerline{\Large{\textbf{Supplementary Materials}}}

\renewcommand{\thesubsection}{\Alph{subsection}}
\section*{Appendix}

\subsection{Notation}
\label{app:notation}

Matrices, vectors and scalars are denoted by uppercase bold, lowercase bold and lowercase normal letters, i.e., $\vX$, $\vx$ and $x$, respectively. The set of $D$-dimensional real numbers is defined by $\bR^D$. $\vx\in\bR^{D}$ and $\vz\in\bZ^{D}$ denote vectors of real and integer variables with dimension $D$. $[\vx ; \vz]$ and $[\vx, \vz]$ indicates column and row concatenation, respectively. $\vX^T$ denotes the transpose of matrix $\vX$ and $[\vX]_j$ its $j$-th row. $\vI_D$ and $\vzero_{D}$ denote the $D$-dimensional identity and zero matrix, respectively. The number of variables is denoted as $\Dz$ and the number of constraints as $\Dc$. We use a hat to indicate a variable prediction, e.g.,  $\hat{z}$. We use the subscript $j$ as a variable index and $i$ as a constraint index. We use uppercase calligraphic and Greek letters to define sets. We denote by $\cS$ a set of time series (indexed by $s$) and by $\cT_s = \{1, \ldots, |\cT_s|\}$ the set of timesteps of time series $s$ (indexed by $t$). $\phi$ represents the set of parameters of an optimization problem and $\Theta$ the set of parameters of a NN. The symbol $\bone_{\text{condition}}$ denotes the indicator function with value 1 if condition is true and 0 otherwise. $\KL\left(p || q \right)$ denotes the KL-divergence of distributions $p$ and $q$.

\subsection{Parameter Normalization \& Distribution Shift}
\label{app:distr_shift}
As mentioned in Section~\ref{sec:model}, the parameter normalization (see~Eq. \eqref{eq:param_norm}) makes our model scale invariant, i.e., our trained model can be applied to new instances where the magnitudes of the constraint matrix and cost vector are entirely different from what was seen during training. More importantly, our trained model can also be used in solving MIP problems whose sizes (i.e., length of the cost vector $D_z$) are different from that of training instances. However, during test time, if the problem size grows (or shrinks), the parameters will become smaller (or larger), causing a distribution shift. 

The distribution shift is caused since in training time with the normalization given in~\eqref{eq:param_norm} (with $p=2$) the parameters $a_{ij}$ and $b_i$ are proportional to $\frac{1}{\sqrt{D_z + 1}}$, and $c_j$ is proportional to $\frac{1}{\sqrt{D_z}}$, 
whereas at test time these would be proportional to $\frac{1}{\sqrt{\bar{D}_z + 1}}$ and $\frac{1}{\sqrt{\bar{D}_z}}$, respectively.
This can be rectified by a further rescaling of parameters during test time.

More precisely, given a model that is trained on instances of size $D_z$ and applied on test instances of size $\bar{D}_z$, we normalize the test instance's parameters as

\begin{equation}\label{eq:param_norm_generic}
a_{i,j}\! =\!  \sqrt{\frac{\bar{D}_z + 1}{D_z + 1}} \frac{a_{i, j}}{\| [\va_i^T ; b_i] \|_2},  \quad b_i \! = \sqrt{\frac{\bar{D}_z + 1}{D_z + 1}} \! \frac{b_i}{\| [\va_i^T ; b_i]\|_2}, \quad c_j \!= \! \sqrt{\frac{\bar{D}_z}{D_z}} \frac{c_j}{\|\vc\|_2}.
\end{equation}

This reduces to~Eq. \eqref{eq:param_norm} (with $p=2$), if $\bar{D}_z = D_z$.  
With this normalization, the parameters $(a_{i,j}, b_i)$ and $c_j$ during test time are also proportional to $\frac{1}{\sqrt{D_z + 1}}$ and $\frac{1}{\sqrt{D_z}}$ respectively. 

Given this normalization, our model can also be trained on MIP instances of different sizes. We just need to fix one size $D_z$ as the reference size and for every instance of different size $\bar{D}_z$, the parameters would have a further rescaling as  in~\eqref{eq:param_norm_generic}.

\subsection{Graph Features}
\label{app:feature}

An integral part of \ours\ is the creation of the features of each node that form the matrix $\vU \in\bR^{|\cV| \times D_u}$, which is fed into the GCN (see Eq.\ \eqref{eq:GNN_embedding}), where $D_u$ is the feature dimension. We differentiate the feature construction between variable and constraint nodes.

The set of parameters that the $j$-th variable is associated with are: (i) its cost in the objective $c_j$, (ii) all its parameters in the constraint matrix $\vA$, i.e.,  $\va_j = [\vA^T]_j$ and (iii) the right-hand side of the constraints $\vb$. 
Usually the constraint matrix $\vA$ is sparse, especially for large problems, which means that each variable participates in only a small subset of the constraints. We leverage this observation and we include as relevant features only the subset of constraints that each variable is included, i.e., for the $j$-th variable we include the constraint $i$ if $a_{i,j} \neq 0$. In particular, we calculate the maximum number of constraints that any variable appears, $m_c$, and we keep that many constraint parameters for all variables (so all the feature dimensions are consistent).\footnote{The maximum number of constraints per variable $m_c$ is calculated across the whole dataset. This number can vary in different datasets but this does not affect the model since no operation is done across this dimension.}

Now, for the $j$-th variable we construct $m_c$ triplets $(a_{i, j}, b_i, c_j)$, where $i \in \cI_j = \{i | a_{i,j} \neq 0\}$. If a variable is included in less that $m_c$ constraints, i.e., if $|\cI_j|< m_c$, then we include $m_c - |\cI_j|$ zero triplets. Each triplet is mapped to the target feature dimension $D_u$ using an MLP with parameters learned jointly across all variables:
\begin{equation}
	\vV_{j} = f^{(\text{var-map})}(\va_{\cI_j, j}, \vb_{\cI_j}, c_j) \!= \! \begin{bmatrix}
		a_{i_1, j} & b_{i_1} & c_j  \\
		\vdots & \vdots & \vdots\\
		a_{i_{m_c}, j} & b_{i_{m_c}} & c_j 
	\end{bmatrix} \vW^{(\text{var-map})},\notag
\end{equation} 
where $\vW^{(\text{var-map})} \in \bR^{3 \times D_u}$ is the learnable weight matrix of the MLP and $\vV_{j} \in \bR^{m_c \times D_u}$. Finally, we aggregate over all rows using a symmetric function (e.g., average):
\begin{equation}
	\vu_{j}^{(\text{var})} =  f^{(\text{var-agg})}(\vV_{j}),
\end{equation} 
where $\vu_j^{(\text{var})} \in \bR^{D_u}$. Note that $\vu_j^{(\text{var})} $ is variable and constraint permutation invariant. The complete feature matrix for the variable nodes is constructed by concatenating all the feature vectors $\vu_j^{(\text{var})}$, i.e., 
$\vU^{(\text{var})} = [\vu_1^{(\text{var})}, \dots, \vu_{D_z}^{(\text{var})}] \in \bR^{\Dz \times D_u}$.

The features of the constraints are constructed in a similar manner. The set of parameters that the $i$-th constraint is associated with are: (i) all its parameters in the constraint matrix $\vA$, i.e.,  $\va_i = [\vA]_i$ (ii) the right-hand side of the constraint $b_i$ and (iii) the costs in the objective $\vc$ for the variables in the constraint. As in the variable nodes, we leverage the sparsity of $\vA$ and find the maximum number of variables in all the constraints, $m_v$. Now, for the $i$-th constraint we create $m_v$ triplets $(a_{i,j}, b_i, c_j)$, where $j\in\cJ_i = \{j: a_{i, j} \neq 0\}$. The triplets are projected to the target dimension $D_u$ and then aggregated in the same way as the variable features: 
\begin{gather}
	\vC_{i} = f^{(\text{con-map})}(\va_{i, \cJ_i}, b_i, \vc_{\cJ_i}) \! = \!\! \begin{bmatrix}
		a_{i, j_1} & b_i & c_{j_1}\\
		\vdots & \vdots & \vdots \\
		a_{i, j_{m_v}} & b_i & c_{j_{m_v}}
	\end{bmatrix} \vW^{(\text{con-map})}, \notag \\
	\vu_{i}^{(\text{con})} =  f^{(\text{con-agg})}(\vC_{i}), \notag
\end{gather} 
where $\vW^{(\text{con-map})} \in \bR^{3 \times D_u}$ is the learnable weight matrix of the MLP, $\vC_{i} \in \bR^{m_v \times D_u}$ and $\vu_i^{(\text{con})} \in \bR^{D_u}$. Note that $\vu_i^{(\text{con})}$ is variable and constraint permutation invariant. The complete feature matrix for the constraint nodes is constructed by concatenating all the feature vectors $\vu_i^{(\text{con})}$, i.e., 
$\vU^{(\text{con})} = [\vu_1^{(\text{con})}, \dots, \vu_{D_c}^{(\text{con})}] \in \bR^{\Dc \times D_u}$.

It is possible to differentiate between the equality and inequality constraints and treat them separately. In that case, a different set of triplets is constructed for each constraint and a different MLP can be used to map the triplets to the feature dimension.

The final feature matrix $\vU$ is:
\begin{equation}
	\vU = \begin{bmatrix}
		\vU^{(\text{var})}  \\
		\vU^{(\text{con})} 
	\end{bmatrix} \in\bR^{|\cV| \times D_u}. 
\end{equation}

\subsection{MLP for Continuous Variables}
\label{app:continuous}
Given an assignment of integer variables, we can re-write the MILP~\eqref{eq:binary_MILP} as the following parametric linear program:
\begin{equation}
	\begin{aligned} 
		V(\vzb) = \underset{\vz}{\text{minimize}} & \quad (\vcb)^T \vzb + (\vcc)^T \vzc\\
		\text{subject to} & \quad \vAb \vzc \leq  \vb - \vAb \vzb. \\
	\end{aligned}
	\label{eq:parametric_LP} 
\end{equation}
where $\vz = [\vzb ; \vzc]$, $\vA = [\vAb, \vAc]$, and  $\vc = [\vcb, \vcc]$. The above function $V : \mathbb{R}^{D_z} \rightarrow \mathbb{R}$ is a piecewise convex function~\cite{MPLP, borrelli2017predictive}, and therefore it can be approximated by an MLP with ReLU activation units. This fact motivated us to use an MLP to map the embedding from the GCN to the continuous variables. Notice that also a convex layer~\cite{agrawal2019differentiable} can be used to approximate optimization programs. However, training convex layers is expensive as it requires to differentiate through the optimization program and in our testing it significantly increased the time needed for training. For this reason, we opted to use an MLP instead.


\subsection{On the Bernoulli and Beta Distributions}
\label{app:bernoulli}

The goal of our model is to learn the value of each binary variable as well as measure how
reliable is each prediction. A natural approach to assess how reliable or risky is a decision is via the variance. Recall that a Bernoulli distribution, defined as $\text{Ber}(\pi)$, has a mean $\pi$ and variance $\pi(1 - \pi)$. Thus, for a given mean, the variance is completely defined and fixed, i.e., we cannot have different variance values for a given mean value. This essentially renders the variance uninformative (given the mean) for the purpose of assessing how reliable is a binary variable prediction. On the other hand, in the Beta distribution a fixed mean does not correspond to a fixed variance value and therefore can be used (along with the mean) as a measure to select trustworthy variables. For example, when selecting the $\rho\%$ most reliable variable predictions, one might prefer to include a variable with mean 0.2 and small variance (and fix it to the integer value 0) than a variable with mean 0.1 and very large variance, which indicates that the model has low confidence on that particular prediction although its mean is closer to 0.

\subsection{Training}

\subsubsection{Clenshaw-Curtis Quadrature}
\label{app:ccq}

The Clenshaw-Curtis quadrature \cite{clenshaw1960method} is a method for numerical integration that is based on the expansion of the integrand in terms of Chebyshev polynomials, where the function $f(x)$ to be integrated over the fixed interval $[-1,1]$ is evaluated at the $K$ roots of a Chebyshev polynomial. Then, the integral can be approximated as:
\begin{equation}
	\int_{-1}^{1}f(x)dx \approx \mathbf{w}^\top\mathbf{y}, 
\end{equation}
where 
\begin{equation}\label{eq:quadrature_weights}
	\mathbf{w} = \mathbf{D}^\top\mathbf{d}, 
\end{equation}
are the quadrature weights, with  $\mathbf{d}\in\mathbb{R}^{K/2+1}$ and $\mathbf{D}\in\mathbb{R}^{(K/2+1)\times(K/2+1)}$. In particular, the $k$-th element of $\mathbf{d}$ (with zero-based indexing) is given by
\begin{equation}
	d_k = \begin{cases}
		1, \quad & k=0,\\
		2/(1-(2k)^2), \quad & k=1,\dots,K/2-1,\\  
		1/(1-K^2), \quad & k=K/2,
	\end{cases} 
\end{equation}
while the $(m,k)$-th element of $\mathbf{D}$ is defined as
\begin{equation}
	D_{mk} = \frac{2}{K}\cos\left(\frac{mk\pi}{K/2}\right)\times
	\begin{cases}
		1/2, \quad & k=0, K/2,\\
		1, \quad & \text{otherwise}.
	\end{cases}
\end{equation}
Notice that both $\mathbf{D}$ and $\mathbf{d}$ are independent of the function $f(x)$ and therefore can be precomputed in $O(K\log K)$. The information of the function is encoded in the vector  $\mathbf{y}\in\mathbb{R}^{K/2+1}$, where its $k$-th element, with $k=0, \dots, K/2$, can be computed as
\begin{equation}\label{eq:y_k}
	y_k = f(\cos(k\pi/K)) + f(-\cos(k\pi/K)).
\end{equation}

Finally, if the function needs to be integrated in the $[0,1]$ interval, we can apply the simple change of variable $x^\prime= 2x-1 \rightarrow x = (x^\prime + 1)/2$ and the new integral becomes
\begin{equation}\label{eq:new_var}
	\int_{-1}^1 \frac{1}{2}f((x^\prime + 1)/2) d\x^\prime,
\end{equation}
where we have also accounted for the Jacobian factor which is equal to $1/2$.  The integral in \eqref{eq:new_var} has the right limits and the Clenshaw-Curtis quadrature can be applied accordingly.

\subsubsection{Regularization}
\label{app:reg}

Ideally, we would like the model to produce Beta distributions with mean close to 0.5 and high variance if it is not confident. We propose the following regularization term to reinforce this behaviour:
\begin{equation}
\label{eq:reg}
	\begin{aligned}
		r(\alpha, \beta, z^\star) = & \int_{0}^{1} \Beta(\pi ; \alpha, \beta) |z^\star - \pi | d\pi \cdot  \KL\left(\text{U}(0, 1) || \Beta(\pi, \alpha, \beta)\right) \\
		= & \left((1 - z^\star) \int_{0}^{1} \Beta(\pi ; \alpha, \beta) \pi d\pi + z^\star \int_{0}^{1} \Beta(\pi ; \alpha, \beta) (1 - \pi) d\pi \right) \\
		& \times  \int_{0}^{1} \log\left(\frac{1}{\Beta(\pi ; \alpha, \beta)}\right)d\pi \\
		= & \left((1 - z^\star) \frac{\alpha}{\alpha + \beta}  + z^\star \frac{\beta}{\alpha + \beta}  \right) \left(- \int_{0}^{1}\log  \left( \frac{\pi^{\alpha - 1} (1-\pi)^{\beta - 1} }{B(\alpha, \beta)} \right) d\pi  \right) \\
		= & \left( \frac{(1 - z^\star) \alpha + z^\star \beta}{\alpha + \beta} \right) \left( \alpha -1 + \beta -1 + \log \left(B(\alpha, \beta) \right)  \right),
	\end{aligned}
\end{equation}
where $\text{Uni}(0, 1)$ is the uniform distribution in the unit interval and $B(\alpha, \beta)$ the Beta function. The goal is to minimize the KL-divergence between the learned Beta distribution and the Uniform 
only when a prediction is wrong. The term  $\int_{0}^{1} \Beta(\pi ; \alpha, \beta) |z^\star - \pi| d\pi$ evaluates how far is the distribution from the true label $z^\star$: if Beta \emph{is close} to the label and has small variance (confident and correct) then this term is very small and the regularization becomes negligible. If the distribution \emph{is far} from the label, the KL divergence regularizes Beta towards Uniform (pushes its mean towards 0.5) and increases its variance.

The overall supervised loss is given by:
\begin{equation}
\ell_{\text{sup-reg}} = \ell_{\text{sup}}(\vz^\star, \psi=\{\valpha_{s,t}, \vbeta_{s,t} \})  + \lambda_{\text{reg}}r(\vz^\star, \psi=\{\valpha_{s,t}, \vbeta_{s,t} \}), \label{eq:supregloss}
\end{equation}
where $\lambda_{\text{reg}} \geq 0$ is a regularization parameter. 


\subsubsection{Weighted Loss}
\label{app:weight}

It is common in classification tasks to include weights in the loss when the classes are imbalanced or to trade off precision and recall. 

Our problem can be interpreted as a classification task with $2^{D_z^{(\text{b})}}$ classes, one for each possible binary sequence of length $D_z^{(\text{b})}$, and each instance then belongs to a specific class. With this interpretation the number of classes explodes very fast and no meaningful technique can be applied since the samples required to have a decent representation of each class in unattainable (e.g., for $D_z^{(\text{b})}=100$ and assuming an average of only 10 samples per class we would require $2^{100} \cdot 10 \approx 10^7$ samples). 

An interpretation that scales is to consider our problem as $D_z^{(\text{b})}$ classification tasks (one fore each variable) with 2 classes each (since we consider binary variables). With this approach, each instance is a combination of $D_z^{(\text{b})}$ classes. Now, each classification task has as many samples as the number of instances, which are split in a different way among the 2 classes.

Note that the above interpretations do not change the mathematical formulation of the problem or the loss but the way we think of the class weighting. Following the second approach, we can calculate the class representation percentage $r_j$ of the $j$-th variable (i.e., each classification task) from the available instances (labels). For example, if in $80\%$ of the labels $z_j=0$ and in $20\%$ of the labels $z_j=1$, then $r_j=0.2$. Then, the (supervised) loss in its generic form would be:
\begin{equation}
	\ell(\vz^\star, \phi) = -\sum_{s=1}^{S}\sum_{t=1}^{T} \sum_{j=1}^{D_z} \frac{ \log(\cdot)}{r_j^{z_{s,t,j}^\star} (1 - r_j)^{1 - z_{s,t,j}^\star}}. 
\end{equation}

The use (or not) of these weights for the computation of the supervised loss is a hyperparameter that is optimized for each dataset.

\subsection{Datasets}
\label{app:datasets}

\input{texfiles/datasets-appendix.tex}

\subsection{Experimental Details}
\label{app:exp_details}

\subsubsection{Benchmarks}
In our experiments we use as benchmark the \nd\ model \cite{nair2020solving} as it is the most related approach to our method \ours.

\nd\ uses a GCN based on the bipartite graph representation of the MILP, similar to our approach, while the output of the GCN is mapped to Bernoulli parameters using an MLP. The model is trained by minimizing the NLL (cross-entropy) modified based on the SelectiveNet approach \cite{geifman2019selectivenet}. The idea is to learn two sets of binary variables: one set indicates if a variable is going to be selected and the other its binary value (to be used if selected). Therefore, the overall loss can be seen as a weighted cross-entropy. Further, this loss assumes a predefined target selection percentage and therefore a different model is required for different variable selection percentages. 

We implemented \nd\ by keeping the common structure of our model and changing the final layers and the loss. Further, we used the same set of features (see Appendix \ref{app:feature}) as in our method for a fair comparison. Note that in \cite{nair2020solving} it is not clearly mentioned what features are used for the \nd\ method. The authors only mention a set of solver-based features  that can be used in their \texttt{neural-branching} method but are not applicable for \nd\ since no such features are available during inference. 

For the method from~\cite{bertsimas2021voice}, we use the implementation available online at \texttt{\hyperlink{github.com/bstellato/mlopt}{github.com/bstellato/mlopt}}. This strategy learns a mapping from problem parameters to feasible variable assignments. First, the input data is categorized into strategies, i.e., feasible variables assignment. Then, the mapping is learned as a multi-class classification problem. Notice that this approach does not generalize when a class is not in the training dataset, i.e., when at test time all feasible assignments seen during training are note feasible for the new problem instance. In order to scale to large datasets, we modified their implementation so that it accepts sparse training data. Moreover, we made sure that their method receives the same train/validation/test split as our method.

\subsubsection{Hyperparameters}

For all datasets we used the Adam optimizer for the gradient updates with $10^{-5}$ weight-decay and clipped the gradient norm to 10. The learning rate was warmed-up linearly from $1\times10^{-4}$ to $1\times10^{-2}$ for the first $500$ steps after which a cosine decay follows for the remaining time steps with a decay rate of 0.99.

\begin{table*}
	\small
	\center
	\caption{\small Hyperparameters of the loss components.}
	{\begin{tabular}{c| ccc}
			\toprule
			Schedule   &        \multicolumn{3}{c}{Parameter}    \\
			\midrule
			& $\lambda$ &  $\lambda_{\text{reg}}$  &  $\lambda_{\text{c}}$   \\
			\midrule
			\midrule
			warm-up steps & \{500, 1000, 1500, 2000\} & $\{250, 500\}$ & $\{1000, 1500\}$   \\
			warm-up initial value & $\{0.01, 0.1\}$ & $\{0.01, 0.1\}$ & $0.1$  \\
			warm-up final value & $\{0.1, 1\}$ & $\{0.1, 1\}$  & $\{1, 10\}$ \\
			final value & $\{1, 10, 50\}$ & $\{0.1, 1, 10, 100\}$  & $\{10, 100\}$ \\
			\bottomrule
	\end{tabular}}
	\label{tab:hp_loss}
\end{table*}

\begin{table*}
	\small
	\center
	\caption{\small Hyperparameters of the model blocks.}
	{\begin{tabular}{c| ccccc}
			\toprule
			Parameter   &        \multicolumn{5}{c}{Block}    \\
			\midrule
			& MLP (features) &  GCN  &  LSTM &   MLP (binary map)  & MLP (continuous map)  \\
			\midrule
			\midrule
			\# of layers & 1 & $\{2, 3\}$ & $\{1, 2, 3\}$ & 1 & 1  \\
			\# of units  & $\{8, 16\}$ & $\{8, 16\}$ & $\{16, 32\}$ & $2\times \Dz^{(\text{b})}$ & $\Dz^{(\text{c})}$  \\
			nonlinearity & $\{\text{ReLU}, \text{No}\}$ & $\{\text{ReLU}, \text{No}\}$  & No & ReLU & No \\
			\bottomrule
	\end{tabular}}
	\label{tab:hp_model}
\end{table*}

For each dataset we performed a hyperparameter optimization on the model and loss parameters. The batch size and the number of training steps were tuned per dataset with values in $\{8, 12, 16, 32\}$ and $\{8000, 15000, 20000\}$, respectively. 
All the loss hyperparameters were linearly warmed-up to an initial value and then linearly scaled up to a final target value in the final training step. The loss scheduling hyperparameter ranges are provided in Table \ref{tab:hp_loss}.
The model hyperparameter ranges are provided in Table \ref{tab:hp_model}.\footnote{Note that the number of units in GCN increases due to skip connections.}

\subsubsection{Computational Environment}
For running our experiments we used p3.8xlarge AWS EC2 instances with 4
Tesla V100 GPUs, 36 CPUs, and 244 GB of memory, and m5.8xlarge	instances with 32 CPUs and	128	 GB of memory.

\subsection{Additional Results}
\label{app:add_exp}

\subsubsection{Additional Datasets}
\label{app:extra_data}

In Table~\ref{tab:add_data} we provide results for accuracy, infeasibility, and optimality gap for the two additional datasets \prs{} and \rems{}, which are smaller versions of \pr{} and \rems{}, respectively. Table~\ref{tab:bertsimassm} shows the accuracy, infeasibility, and optimality gap on these two additional datasets for the method from~\cite{bertsimas2021voice}. We notice that our method performs better than \nd{} and the method from~\cite{bertsimas2021voice} in almost all additional experiments. Notice that the method from~\cite{bertsimas2021voice} is able to find a feasible solution for \prs{} in more cases than \ours{} and \nd{}. However, this method is always infeasible in \rem{}, \ours{} and \nd{} always find a feasible solution.

\begin{table*}[t]
	\center
	\caption{\small Accuracy, infeasibility and optimality gap (mean \textpm\ std) of \ours{} vs. \nd{}. Bold indicates best method.}
	\scriptsize
	{\begin{tabular}{cc | cc |cc | cc}
			\toprule
			 & & \multicolumn{2}{|c}{\textbf{Accuracy}} & \multicolumn{2}{|c}{\textbf{Infeasibility}} & \multicolumn{2}{|c}{\textbf{Optimality gap}} \\
			\midrule
			$\rho$  & Method   &       \prs &          \rems &           \prs  &        \rems &        \prs &       \rems   \\
			\midrule
			\midrule
			\multirow{3}{*}{30\%} & \ours &  \textbf{99.18 \textpm\ 0.17} &  \textbf{100.00 \textpm\ 0.00} &  \textbf{7.88 \textpm\ 0.63} &  \textbf{0.00 \textpm\ 0.00} &   \textbf{0.35 \textpm\ 0.01} &  \textbf{0.12 \textpm\ 0.03} \\
			& \nd &  95.67 \textpm\ 0.48 &  \textbf{100.00 \textpm\ 0.00} &  77.91 \textpm\ 9.73 &  \textbf{0.00 \textpm\ 0.00} &   \textbf{0.39 \textpm\ 0.04} &  0.16 \textpm\ 0.05  \\
			
			\midrule
			
			\multirow{3}{*}{40\%} & \ours &  \textbf{98.29 \textpm\ 0.22} &  \textbf{100.00 \textpm\ 0.00} &  \textbf{35.63 \textpm\ 21.02} &  \textbf{0.00 \textpm\ 0.00} &   0.95 \textpm\ 0.37 &  \textbf{0.07 \textpm\ 0.02}\\
			& \nd &  78.33 \textpm\ 25.14 &  \textbf{100.00 \textpm\ 0.00} &  96.27 \textpm\ 4.32 &  \textbf{0.00 \textpm\ 0.00} &  \textbf{0.27 \textpm\ 0.02} &  0.11 \textpm\ 0.07\\
			\midrule
			
			\multirow{3}{*}{50\%} & \ours &  \textbf{97.32 \textpm\ 0.40} &  \textbf{100.00 \textpm\ 0.00} &  \textbf{65.07 \textpm\ 29.85} &  \textbf{0.00 \textpm\ 0.00} &  \textbf{1.48 \textpm\ 0.89} &  0.05 \textpm\ 0.02\\
			& \nd &  89.46 \textpm\ 0.93 &  \textbf{100.00 \textpm\ 0.00} &  100.00 \textpm\ 0.00 &  \textbf{0.00 \textpm\ 0.00} &  - &  \textbf{0.03 \textpm\ 0.03}\\
			\midrule
			
			\multirow{3}{*}{60\%} & \ours &  \textbf{95.19 \textpm\ 0.42} &  \textbf{100.00 \textpm\ 0.00} &  \textbf{100.00 \textpm\ 0.00} &  \textbf{0.00 \textpm\ 0.00} &  - &  0.03 \textpm\ 0.03\\
			& \nd &  92.41 \textpm\ 1.26 &  \textbf{100.00 \textpm\ 0.00} &  \textbf{100.00 \textpm\ 0.00} &  \textbf{0.00 \textpm\ 0.00} &  - &  \textbf{0.02 \textpm\ 0.01}\\
			\midrule
			
			\multirow{3}{*}{70\%} & \ours &  \textbf{92.79 \textpm\ 0.43} &  \textbf{100.00 \textpm\ 0.00} &  \textbf{100.00 \textpm\ 0.00} & \textbf{0.00 \textpm\ 0.00} &  - &  \textbf{0.03 \textpm\ 0.04}\\
			& \nd &  89.52 \textpm 2.55 &  \textbf{100.00 \textpm\ 0.00} &  \textbf{100.00 \textpm\ 0.00} &  \textbf{0.00 \textpm\ 0.00} &  - &  \textbf{0.03 \textpm\ 0.01}\\
			\bottomrule
		\end{tabular}
	}
	\label{tab:add_data}
\end{table*}

\begin{table}[t]
	\center
	\caption{\small Accuracy, infeasibility and optimality gap of \cite{bertsimas2021voice}.}
	\scriptsize
	{\begin{tabular}{c | c |c | c}
			\toprule
			 & \textbf{Accuracy} & \textbf{Infeasibility} & \textbf{Optimality gap} \\
			\midrule
            \prs & 85.4 & 93.0 & 9.4 \\
            \rems & 99.22 & 100.00 & $\infty$ \\
			\bottomrule
		\end{tabular}
	}
	\label{tab:bertsimassm}
\end{table}

\subsubsection{Unsupervised Loss}
\label{app:unsup_results}

As we already showed in Section \ref{sec:exp_acc}, training \ours{} only with supervised loss (\ourssup{}) leads to decreased performance for reasons we discussed throughout the paper. Here, we complete the picture by including the performance of \oursuns{}, i.e., \ours{} trained only with unsupervised loss. Tables~\ref{tab:accuracy_uns}-\ref{tab:gap_uns} show the accuracy, infeasibility and optimality gap, respectively. It is clear that the performance of \oursuns{} compared to \ours{} is significantly degraded in most datasets, with the sole exception of the \ca. 

It is important to highlight some key observations based on these results. 
First, the accuracy of \oursuns{} in the \eg\ dataset is less than 50\% for all values of $\rho$. This does not lead to significant infeasibilities since apparently the constraints of the problem are not that hard to satisfy, however it causes a significant optimality gap. This is a case where labels make a big difference in the model training process and in fact the performance of \ourssup{} is in par with \ours{}, i.e., most of the learning is done from the labels. On the other end, we observe that in the \ca\ dataset the model is able to learn the optimal values without labels. Arguably this is an easy problem to learn since all the models achieve optimal performance. In such cases the overhead of collecting (many) labels might not be necessary. All the other datasets fall in the case where both \ourssup{} and \oursuns{} underperform compared to \ours{} and a combination of supervised and unsupervised loss gives the best performance.

\subsubsection{Missing Labels}
\label{sec:exp_miss}
\begin{figure*}[t]
	\centering
	\includegraphics[width=0.9\textwidth]{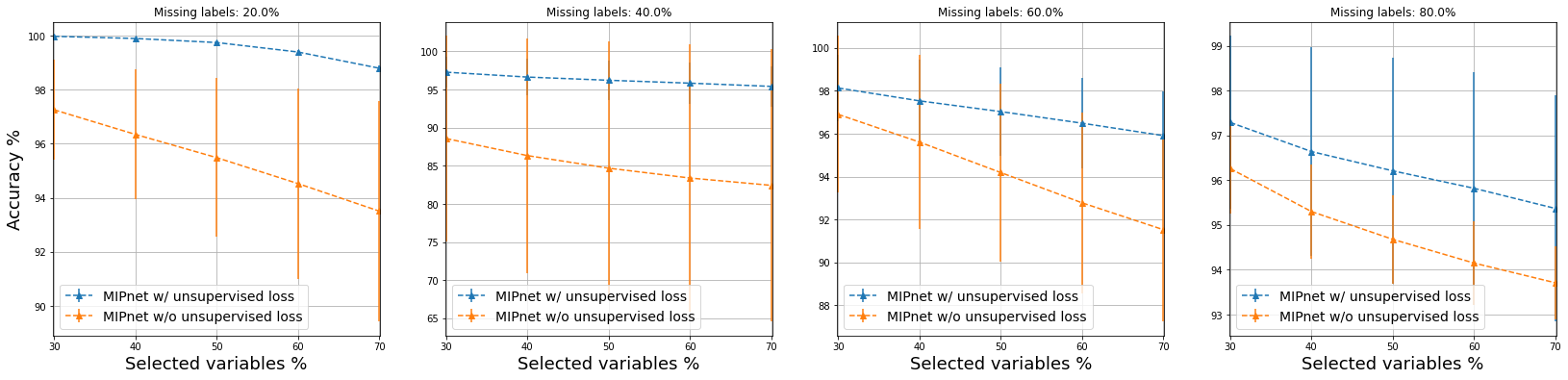}
	\caption{Comparison of \ours\ w\textbackslash\ and w\textbackslash o unsupervised loss for different missing label percentages.}
	\label{fig:miss}
\end{figure*}

Finally, we examine how missing labels affect the performance of a model trained with and without unsupervised loss. For this experiment we used the \tsp{} dataset. Figure \ref{fig:miss} illustrates the accuracy of the models for different percentages of missing labels. In the presence of missing labels the model that is trained with both supervised and unsupervised loss clearly outperforms the one without unsupervised loss in all cases.

\begin{table*}[t]
	\footnotesize
	\center
	\caption{\small Accuracy (mean \textpm\ std)$\%$ of \oursuns{}.}
	\resizebox{0.9\textwidth}{!}{\begin{tabular}{c | cccccc}
		\toprule
		\midrule
		   $\rho \ (\%)$     &       \pr &          \fl &           \tsp  &        \eg &        \rem   &  \ca\\
		\midrule
		\midrule
		30\% &   95.46 \textpm\ 1.47 &  99.02 \textpm\ 1.12 &  96.45 \textpm\ 2.69 &  49.38 \textpm\ 1.36 &  100.00 \textpm\ 0.00 & 100.00 \textpm\ 0.00\\

		\midrule
		
		40\% &   94.13 \textpm\ 1.16 &  98.90 \textpm\ 1.25  &  95.83 \textpm\ 2.90 &  48.16 \textpm\ 0.75 &  100.00 \textpm\ 0.00 & 100.00 \textpm\ 0.00\\
		
		\midrule

		50\% &   92.60 \textpm\ 1.13 &  98.85 \textpm\ 1.26 &  95.28 \textpm\ 3.04 &  47.61\textpm\ 0.77 &  99.99 \textpm\ 0.00 & 100.00 \textpm\ 0.00\\
		
		\midrule
		
		60\% &   91.24 \textpm\ 1.01 &  98.76 \textpm\ 1.21 &  94.76 \textpm\ 3.03 &  47.59 \textpm\ 0.86 &  99.99 \textpm\ 0.00 & 100.00 \textpm\ 0.00\\
	
		\midrule

		70\% &   89.49 \textpm\ 0.76 &  98.50 \textpm\ 1.02  &  94.32 \textpm\ 2.84 &  47.78 \textpm\ 0.79 &  99.99 \textpm\ 0.00 & 100.00 \textpm\ 0.00\\
	
		\bottomrule
	\end{tabular}}
	\label{tab:accuracy_uns}
\end{table*}

\begin{table*}[t!]
	\footnotesize
	\center
	\caption{\small Infeasibility (mean \textpm\ std)$\%$ of \oursuns{}.}
	\resizebox{0.9\textwidth}{!}{\begin{tabular}{c | cccccc}
		\toprule
		\midrule
		   $\rho \ (\%)$     &       \pr &          \fl &           \tsp  &        \eg &        \rem   &  \ca\\
		\midrule
		\midrule
		30\% &   100.00 \textpm\ 0.00 &  6.92 \textpm\ 0.56 &  33.43 \textpm\ 47.07 &  0.02 \textpm\ 0.03 &  0.00 \textpm\ 0.00 & 0.00 \textpm\ 0.00\\

		\midrule
		
		40\% &   100.00 \textpm\ 0.00 &  13.87 \textpm\ 10.09  & 41.94 \textpm\ 42.22 &  0.06 \textpm\ 0.00 &  0.00 \textpm\ 0.00 & 0.00 \textpm\ 0.00 \\
		
		\midrule

		50\% &   100.00 \textpm\ 0.00 &  29.89 \textpm\ 25.36 &  54.68 \textpm\ 39.66 &  0.11 \textpm\ 0.06 &  0.00 \textpm\ 0.00 & 0.00 \textpm\ 0.00 \\
		
		\midrule
		
		60\% &   100.00 \textpm\ 0.00 &  47.28 \textpm\ 24.40 &  69.16 \textpm\ 40.00  & 0.15 \textpm\ 0.03 &  0.00 \textpm\ 0.00 & 0.00 \textpm\ 0.00 \\
	
		\midrule

		70\% &   100.00 \textpm\ 0.00 &  72.68 \textpm\ 15.37  &  84.04 \textpm\ 21.92 &  0.17 \textpm\ 0.08 &  0.00 \textpm\ 0.00 & 0.00 \textpm\ 0.00 \\
	
		\bottomrule
	\end{tabular}}
	\label{tab:infeas_uns}
\end{table*}

\begin{table*}[t!]
	\footnotesize
	\center
	\caption{\small Optimality gap (mean \textpm\ std)$\%$ of \oursuns{}.}
	\resizebox{0.9\textwidth}{!}{\begin{tabular}{c | cccccc}
		\toprule
		\midrule
		   $\rho \ (\%)$     &       \pr &          \fl &           \tsp  &        \eg &        \rem   &  \ca\\
		\midrule
		\midrule
		30\% &   - &  8.83 \textpm\ 12.33 &  1.55 \textpm\ 1.27 &  8.59 \textpm\ 0.04 &  0.00 \textpm\ 0.00 & 0.00 \textpm\ 0.00\\
	
		\midrule
		
		40\% &  - &  12.44 \textpm\ 17.31  &  3.11 \textpm\ 2.78 &  12.62 \textpm\ 0.36 &  0.00 \textpm\ 0.00 & 0.00 \textpm\ 0.00\\
		
		\midrule

		50\% & - &  15.08 \textpm\ 20.06 &  5.95 \textpm\ 5.38 &  17.16 \textpm\ 0.63 &  0.66 \textpm\ 0.01 & 0.00 \textpm\ 0.00\\
		
		\midrule
		
		60\% & - &  18.97 \textpm\ 24.36 &  11.06 \textpm\ 9.82 & 22.00 \textpm\ 0.95 & 1.23 \textpm\ 0.01 & 0.00 \textpm\ 0.00\\
	
		\midrule

		70\% &  - &  32.20 \textpm\ 27.12  &  19.73 \textpm\ 17.64 &  27.41 \textpm\ 1.11 &  1.57 \textpm\ 0.07 & 0.00 \textpm\ 0.00\\
	
		\bottomrule
	\end{tabular}}
	\label{tab:gap_uns}
\end{table*}

\subsection{Discussion of Related Work}\label{app:related}


Table \ref{tab:related-work} summarizes the key differences of \ours{} from the most-related models and approaches. In the evaluation section we have presented detailed comparisons of \ours{} with \cite{nair2020solving}, \cite{geifman2019selectivenet} and \cite{bertsimas2021voice}, while here we focus on qualitative differences. In detail, we see that \ours{} is the only model that uses semi-supervised learning. This is very important since obtaining labels requires, in most cases, solving large-scale NP-hard problems, while the unsupervised loss component can improve the performance of the model further, by allowing us to identify the effect (in terms of optimality gap and infeasibility) of mispredictions at a per-variable granularity. These benefits are evident in the experiments. We note that other works, such as \cite{nair2020solving} and \cite{dingaaai20}, identify the importance of unsupervised learning as well, but this approach is esssentially explored only in \cite{karalias202Nips}. 

In particular, \cite{karalias202Nips} proposes an unsupervised approach specifically for graph problems where a GNN is used to learn a distribution over the graph nodes, representing a solution. In order to produce an integral solution, the authors derandomize the continuous values using sequential decoding. The method comes with theoretical guarantees, giving good and feasible solutions with high probability. However, this approach is not applicable to general MIPs, while even for the targeted graph problems it is not trivial to include general constraints. Besides, this work does not benefit from the availability of labels, which, in certain operational environments are available. 

On the other hand, in \cite{nair2020solving} the authors provide only a preliminary evaluation with training datasets that might include noisy labels (obtained from suboptimal solutions); while \cite{dingaaai20} performs an initial variable assignment and then uses neighborhood search to increase the training data without solving the problem exactly -- which works under the assumption of locality, as the authors stress. On the contrary, \ours{} uses a \emph{semi-supervised} learning approach which, apart from allowing us to expand the training dataset, enables tuning the impact of each variable on the feasibility and objective value of the problem.

Another distinguishing feature of \ours{} is that, unlike all prior works, it accounts for the temporal structure across the different instances, which, as shown in the experiments, indeed enhances the model's performance. This aspect is crucial as, more often than not, the different instances that practitioners solve exhibit a temporal structure, namely problem parameters such as the commodity volumes, user requests, transportation costs or electricity prices, follow some type of diurnal pattern.

We also note that the related works in Table \ref{tab:related-work} can be separated to those that use GCN in order to automate feature embedding, and to those that do not follow this approach, as e.g., \cite{bertsimas2019online, bertsimas2021voice}. GCNs seem to enhance the performance and provide the means to account for dependencies across the variables and, indirectly, among the constraints as well. Finally, \ours{} is the only work that uses a Bayesian approach and a tunable confidence for performing the assignment of variables, unlike, e.g., the threshold-based variable selection rule of \cite{dingaaai20}. We note also that \cite{nair2020solving} offers the option to select different percentage $\rho$ of the binary variables, but this decision needs to be made in advance so as to train the model accordingly.

\begin{table*}[h!]
	\scriptsize
	\center
	\caption{\small Comparison of \ours{} with most related models.}
	{\begin{tabular}{c | cccc}
		\toprule
		\midrule
    & \textbf{Uses GCN} &\textbf{Vars Select Method} &\textbf{Temporal} &\textbf{Learning}   \\
		\midrule
		\midrule
    \cite{nair2020solving}   & Yes       & Bernoulli &  No          &  Superv.         \\
		\midrule
    \cite{bertsimas2021voice} & No       & No        &  No          &  Superv.         \\
		\midrule
    \cite{dingaaai20}        & Yes       & Threshold-based &  No          &  Superv.         \\
		\midrule
    \cite{bertsimas2019online}        & No       & No &  No          &  Superv.         \\
		\midrule
    \cite{gasse2019exact}        & Yes       & Threshold-based &  No          &  Superv.         \\
		\midrule
    \cite{karalias202Nips}        & Yes       & No &  No          &  Unsuperv.         \\
		\midrule
    \ours{}        & Yes       & Bayesian (\& tunable thresh.) &  Yes          &  Superv. + Unsuperv.         \\
		\midrule
    \bottomrule
    \bottomrule
	\end{tabular}}
	\label{tab:related-work}
\end{table*}


%% file: texfiles/datasets-appendix.tex

We introduce the six types of problems which we employed in the evaluation of \ours{}. The results for the problems from Tables~\ref{tab:datasets} and \ref{tab:timeseries} can be found in Section~\ref{sec:experiments}, and the results from problems from Tables~\ref{tab:datasets_add} and \ref{tab:timeseries_additional} are discussed in Section~\ref{app:extra_data} of the Appendix. These problems were selected in line with the evaluation studies in prior works, see \cite{nair2020solving, dingaaai20, bertsimas2021voice} and references therein; they appear regularly in different application domains and, furthermore, constitute building blocks for a very wide range of optimization problems. In the sequel we describe the mathematical formulation of each problem, denoted as $\mathbb{P}_{p}^{t}$ with $t\in\cT$ and $p=1, \dots, 6$, to indicate that each formulation corresponds to a particular instance of the problem with time-varying parameters.\footnote{We have ignored the time series index $s$ to avoid notational cluttering.}\footnote{For the problem formulations we have used $t$ as a superscript in all the time-varying parameters ($t$ is fixed per instance and it just indicates that a specific parameter changes across instances), while all the subscripts are instance-dependent indices.} We further identify the problems' complexity and class (type of problem), and we map them to real-world applications. Finally, we provide details about the datasets and instances we used in training, testing and validation -- see also the summarized view in Tables~\ref{tab:datasets}--\ref{tab:timeseries_additional}.

\begin{table*}
	\small
	\center
	\caption{\small Datasets problem sizes.}
	{\begin{tabular}{c | cccccc}
		\toprule
		\midrule
		         &    \pr &          \fl &           \tsp  &        \eg  &        \rem  & \ca  \\
		\midrule
		\midrule
		    \# Integer variables        &       2812  &          1989 &           144  &        1000 &        30000   & 39942\\
		\midrule
		    \# Continuous variables   &       0 &          0 &                  500 &     0 &   0  & 0 \\
		\midrule
		    \# Equality constraints   &       450 &          50 &           36  &        1 &    0   & 0\\
		\midrule
		    \# Inequality constraints  &       36  &          39 &           4082  &        1010 &        20   & 1\\
		\midrule
		    Data                         &       real  &          mixed &           synthetic  &        synthetic & synthetic   & mixed  \\ \midrule
Fraction of non-zeros in $\vA$               &       0.0374  &          0.0226 &           0.1414  &        0.6723 &        0.50  & 1.0 \\ \midrule
		    Fraction of non-zeros in $\vb$              &       1.0        &          0.5       &           0.8333  &        0.9876  &        1.0&         1.0  \\ \midrule
		    Fraction of non-zeros in $\vc$                     &       0.579   &          0.9804 &           0.9996  &        1.0  &        1.0   & 0.0187  \\ \midrule
		    Fraction of non-zeros in $\vz^*$               &       0.0001   &          0.0278 &           0.0833  &        0.7375   & 0.0044 & 0.0006   \\ 
		\bottomrule
	\end{tabular}}
	\label{tab:datasets}
\end{table*}

\begin{table*}[t]
	\small
	\center
	\caption{\small Additional datasets problem sizes.}
	{\begin{tabular}{c | cc}
		\toprule
		\midrule
		        &  \prs  & \rems  \\
		\midrule
		\midrule
		    \# Integer variables & 1044        &   10000 \\
		\midrule
		    \# Continuous variables & 0   & 0 \\
		\midrule
		    \# Equality constraints & 390   & 0\\
		\midrule
		    \# Inequality constraints & 36  &    10\\
		\midrule
		    Data                      & real   & synthetic \\ \midrule
Fraction of non-zeros in $\vA$         & 0.0333      &        0.501  \\ \midrule
		    Fraction of non-zeros in $\vb$   & 1.0            &       1.0 \\ \midrule
		    Fraction of non-zeros in $\vc$         & 0.9169             &       1.0    \\ \midrule
		    Fraction of non-zeros in $\vz^*$   &       0.4298            &      0.0166 \\ 
		\bottomrule
	\end{tabular}}
	\label{tab:datasets_add}
\end{table*}

\begin{table*}[t]
	\small
	\center
	\caption{\small Datasets time series dimensions.}
	{\begin{tabular}{c | cccccc}
		\toprule
		\midrule
		        &       \pr &          \fl &           \tsp  &        \eg &        \rem & \ca  \\
		\midrule
		\midrule
		    \# of time series for training         &       46  &          80 &           80  &        80 &          80 & 80  \\
		\midrule
		    \# of time series for testing    &       5 &          20 &           20  &        20 &        20  & 20 \\
		\midrule
		    \# of time series for validation    &       5 &          20 &           20  &        20 &        20 & 20  \\
		\midrule
		    \# of instances per time series   &       96  &          700 &           400  &        100 &        100 & 150   \\
		\bottomrule
	\end{tabular}}
	\label{tab:timeseries}
\end{table*}

\begin{table*}
	\small
	\center
	\caption{\small Additional datasets time series dimensions.}
	{\begin{tabular}{c | cc}
		\toprule
		\midrule
		        &       \prs &        \rems   \\
		\midrule
		\midrule
		    \# of time series for training         &       96 &  80    \\
		\midrule
		    \# of time series for testing    &       18 &          20  \\
		\midrule
		    \# of time series for validation    &       10 &          20   \\
		\midrule
		    \# of instances per time series   &       96  &          100    \\
		\bottomrule
	\end{tabular}}
	\label{tab:timeseries_additional}
\end{table*}

\subsubsection{Multicommodity Network Design and Routing}

We consider the joint unsplittable-routing and network design problem where the goal is to minimize the aggregate routing cost by deciding how much capacity to purchase for each link and which path to select for each commodity. Formally, the network is modeled by a graph $\cl G=(\cl V,\cl E)$ and serves a set of $\mathcal K$ commodities with time varying demands $\bm d^t=\{d^t_k\}_{k\in\cl K}$, and predetermined sources and destinations. Specifically, our goal is to solve the following \emph{binary problem}:
\begin{equation*}
	\begin{aligned}
	\mathbb{P}_{1}^{t}:\quad \underset{\bm{p}, \bm{q}}{\text{minimize}} &\quad  \sum_{k \in \mathcal{K}} \sum_{l \in \mathcal{P}_{k}} v_{k, l} p_{k, l} d_k^{t} + \sum_{e \in \mathcal{E}} \sum_{i \in \mathcal{I}} u_i q_{e, i} \\
		\text{subject to} & \quad \sum_{l \in \mathcal{P}_{k}} p_{k, l} = 1, &&\forall k \in \mathcal{K},\\
		& \quad \sum_{k \in \mathcal{K}} \sum_{l \in \mathcal{P}_{k}} a_{e, k, l}p_{k, l} d_k^{t} \leq {b}_{e} + \sum_{i \in \mathcal{I}}c_i q_{e, i}, &&\forall e \in \mathcal{E}, \\
		& \quad \bm{p}\in \{0,1\}^{|\mathcal K|},\,\, \bm{q} \in \{0,1\}^{|\mathcal E| |\mathcal I|}. \\
	\end{aligned}
\end{equation*}
Problem $\mathbb{P}_{1}^t$ has the following binary \emph{variables}:
\begin{itemize}
	\item $p_{k, l} \in \{0,1\}$: selection of path $l$ for routing commodity $k\in\mathcal K$,
	\item $q_{e, i} \in \{0,1\}$: installment of $i$-type capacity at edge $e \in \mathcal{E}$,
\end{itemize}
and the following \emph{parameters}:
\begin{itemize}
    \item $\mathcal P_{k}$: set of eligible paths for commodity $k\in \mathcal K$,
	\item $v_{k, l}\in \mathbb R_+$: routing cost (per unit of traffic) of path $l\in\mathcal P_{k}$,
	\item $a_{e, l}^{k} \in \{0,1\}$: parameter indicating that edge $e$ is contained in path $l\in \mathcal P_{k}$,
	\item $b_{e} \in \mathbb{R}_+$: initial capacity installed at edge $e \in \mathcal{E}$,
	\item $u_i\in \mathbb{R}_+$: cost of buying $i$-th capacity installment,
	\item $c_i\in \mathbb{R}_+$: capacity of $i$-th installment,
    \item $d_k^t\in \mathbb{R}_+$: demand of commodity $k$ at time $t$.	
\end{itemize}
\textbf{Complexity}. $\mathbb{P}_{1}^t$ is known to be NP-hard even for a single commodity, i.e., $|\mathcal K|=1$, cf. \cite{routing-complexity}.

 \textbf{Datasets}. The problem parameters, including the graph topology, were taken from the SNDlib database \cite{SNDlib10} that contains real-world traffic matrices and network topologies from different communication (and other) networks (e.g., from backbone ISP networks).\footnote{SNDlib is a rich and well-known database which is regularly used in the evaluation of routing algorithms in communication networks.} In particular, for this problem we used the dataset \emph{geant} dataset, the size of which can be seen in Table \ref{tab:datasets}. The capacity and instalment costs are fixed over time, while the demand vector $\bm d^t$ is being updated every $15$ minutes over a period of $4$ months. We note that we have scaled the demand $\bm d^t$ and edge capacity by a factor of $100$ and $1/40$, respectively, to increase the complexity of the problem.\footnote{This scaling renders capacity purchase necessary in order to fulfill the demands.} The number of time series we used, and the number of instances for each time series can be found in Table \ref{tab:timeseries}, where the instances differ in vector $\bm{d}^t$.

\subsubsection{Facility Location}
In the facility location problem, given a set $\mathcal{I}$ of \emph{facility} locations and a set $\mathcal{J}$ of \emph{clients}, the goal is to decide which facilities to open and how to assign clients to the opened facilities. The objective is to minimize the total opening and assignment cost. More formally, we solve the following \emph{binary problem}:
\begin{equation*}
\begin{aligned}
\mathbb{P}_{2}^t:\quad \underset{\vz}{\text{minimize}} &\quad \sum_{i \in \mathcal{I}}\sum_{j\in\mathcal{J}} c_{j,i}d_j^t z_{j,i} + \sum_{i \in \mathcal{I}} f_i x_i \\
\text{subject to} & \quad \sum_{i \in \mathcal{I}} z_{j,i} = 1, &&\forall j \in \mathcal{J}, \\
& \quad \sum_{j \in \mathcal{J}}z_{j,i} \leq 2 |\mathcal{J}| x_i, &&\forall i\in \mathcal I, \\
 & \quad \vz \in \{0,1\}^{|\mathcal{I}| |\mathcal{J}||\mathcal{I}|},
\end{aligned}
\end{equation*}
where the vector of binary variables $\vz = [z_{11}, \ldots, z_{|\mathcal{I}| |\mathcal{J}|}, x_1, \ldots, x_{|\mathcal{I}|}]$.
In the above problem, we introduced the following \emph{binary variables}:
\begin{itemize}
	\item $z_{j,i} \in \{0,1\}$: assignment of client $j\in\mathcal J$ to facility $i\in\mathcal I$,
	\item $x_{i} \in \{0,1\}$: deployment of facility $i\in\mathcal I$,
\end{itemize}
and the following \emph{parameters}:
\begin{itemize}
	\item $f_i\in \mathbb R_+$: cost for opening facility $i\in\mathcal I$,
	\item $c_{j,i}\in \mathbb R_+$: cost for associating client $j\in\mathcal J$ to facility $i\in\mathcal I$,
    \item $d_j^t\in \mathbb R_+$: demand of client $j\in\mathcal J$ at time $t$.
\end{itemize}

\textbf{Complexity}. $\mathbb{P}_{2}^t$ is the metric uncapacitated facility location problem, which is NP-hard and inapproximable below $1.463$ ratio, with the currently-known algorithm achieving a solution within $1.488$ ratio, see \citeAP{ufl}. 

 \textbf{Datasets}. The graph topology is borrowed from the \emph{germany50} dataset \cite{SNDlib10}, the cost parameters $f_i$ are randomly generated from a uniform distribution for each problem instance and the time-varying demand $\bm{d}^t = \{d^t_j\}_{j\in \mathcal J}$ is defined as:
\begin{equation*}
\bm{d}^{t+1} = \max\left\{0, {\bm{A} {\bm{d}^t}} + a_1 \sin(t/{\tt{period}_1})+ a_2 \sin(t/{\tt{period}_2}) + {\bf{w}}\right\}, 
\end{equation*}
for a Hurwitz matrix $\bf{A} \in \mathbb{R}^{|\mathcal{J}|\times|\mathcal{J}|}$ constructed by randomly sampling the eigenvalues $\lambda_i \sim \mathcal{U}(0.98, 0.999)$ and the associated orthonormal eigenvectors $\mathbf{v_i}$\footnote{The eigenvectors are computed randomly sampling from a standard normal and then using the Gram–Schmidt procedure to compute a set of orthonormalized vectors.} for $i \in \{1, \ldots, |\mathcal{J}|\}$. The random variables $a_1 \sim \mathcal{U}(1, 5)$ and $a_2 \sim \mathcal{U}(2, 10)$, and the random vector ${\bf{w}} \sim \mathcal{N}(\mathbf{0}, \mathbf{\Sigma_1})$, where the vector of zeros
$\mathbf{0} \in \mathbb{R}^{|\mathcal{J}|}$ and the covariance matrix $\mathbf{\Sigma_1}\in \mathbb{R}^{|\mathcal{J}| \times |\mathcal{J}|}$ is a diagonal matrix of ones. In the above definition the $\max$ operator is applied component-wise and guarantees nonnegative demand. Finally, the periodicity of the sine functions are ${\tt{period}_1}=20$ and ${\tt{period}_2}=70$.

%
%
%
%
%
%
%
%
%
%
%

\subsubsection{Travelling Salesman}
We consider the Travelling Salesman Problem (TSP) under the Dantzig–Fulkerson–Johnson formulation \citeAP{dantzig-routing}. In detail, we have a set $\mathcal N$ of $N$ cities, which we need to visit exactly once and return to the starting node. The cities are connected with a graph $G=(\mathcal N, \mathcal E)$ which induces routing (or distance) costs $c_{i,j}, \forall (i,j)\in\mathcal E$. Our goal is to minimize the aggregate (travelling) cost of the route. In detail, the formulated \emph{binary problem} is:
\begin{equation*}
\begin{aligned}
\mathbb{P}_{3}^t:\quad \underset{\vz}{\text{minimize}} &\quad \sum_{i \in \mathcal{N}}\sum_{j\in\mathcal{N}} c_{i,j}^t z_{i,j}\\
\text{subject to} & \quad \sum_{i \in \mathcal{N}} z_{i,j} = 1, &&\forall j \in \mathcal{N}, \\
& \quad \sum_{j \in \mathcal{J}}z_{i,j}= 1, &&\forall i\in \mathcal N, \\
& \quad \sum_{i \in \mathcal{S}}\sum_{j \in \mathcal{S}}z_{i,j}\leq |\mathcal S|-1, &&\forall \mathcal S \subseteq\{1,\ldots,N\}, \, |\mathcal S|\geq 2, \\
 & \quad \vz \in \{0,1\}^{N^2}. \\
\end{aligned}
\end{equation*}

The problem includes the \emph{binary variables}:
\begin{itemize}
	\item $z_{i,j} \in \{0,1\}$: inclusion of link $(i,j)\in \mathcal E$ or not ($x_{ij}=0$) in the route,
\end{itemize}
and \emph{parameters}:
\begin{itemize}
	\item $c_{i,j}^t\in \mathbb R_+$: cost for traversing link $(i,j)$ at time $t$.
\end{itemize}
\textbf{Complexity}. $\mathbb{P}_{3}^t$ is NP-hard to approximate within any polynomial factor for general distance functions, while improved approximation ratios are available for certain restricted cases (metric spaces, etc.). We refer the reader to \citeAP{karlin-tsp22} for an up-to-date discussion on the complexity of TSP.

 \textbf{Datasets}. We created instances with $N=12$ nodes. The location of each node is randomly generated using a uniform distribution and we assumed that the graph is fully connected. The cost vector ${\bm{c}^t}\in \mathbb{R}^n$  is given by the distance between each city pair and it changes at each time step as follows: 
\begin{equation*}
\bm{c}^{t+1} = \max\{0, \bm{c}^t + \bf{w}\}, 
\end{equation*}
where ${\bf{w}} \sim \text{U}(d_{min}, d_{max})^n$ and the $\max$ operator is applied component-wise; where $d_{min}, d_{max}$ are positive constants. 

\subsubsection{Revenue Maximization}
We consider the problem of shipping commodities from a set $\mathcal N$ of $N=|\mathcal{N}|$ source nodes (each commodity corresponds to one node) towards their intended destinations over predetermined (overlapping) paths. Our goal is to maximize the revenue from delivering as many commodities as possible, while satisfying the time-varying capacity of each edge $i\in\mathcal I$ that lies along the path of each commodity $n\in\mathcal N$, i.e., $i \ni n$. Formally, we consider the following \emph{binary} problem:
\begin{equation*}
\begin{aligned}
\mathbb{P}_{4t}: \quad \underset{\vz }{\text{maximize}} &\quad \sum_{n\in\mathcal{N}} c^t_{n} z_{n} \\
\text{subject to} & \quad \sum_{n\in\mathcal{N}} a_{i, n} z_{n} \leq b^t_i, \quad \forall i \in \mathcal{I}, \\
 & \quad \vz \in \{0,1\}^{N}. \\
\end{aligned}
\end{equation*}
The problem includes the \emph{binary variables}:
\begin{itemize}
	\item $z_{n} \in \{0,1\}$: transporting commodity $n$ to its destination via a predetermined path, or not,
\end{itemize}
and \emph{parameters}:
\begin{itemize}
	\item $c_{n}^t\in \mathbb R_+$: revenue when delivering commodity $n$, at time $t$,
	\item $b_{i}^t\in \mathbb R_+$: transportation capacity of link $i\in\mathcal I$, at time $t$,
	\item $a_{i,n} \in \{0,1\}$: indicating if link $i$ lies in the predetermined path of commodity $n$, i.e., $i\ni n$.
\end{itemize}

\begin{figure}[t]
    \centering
    \includegraphics[width=0.5\textwidth]{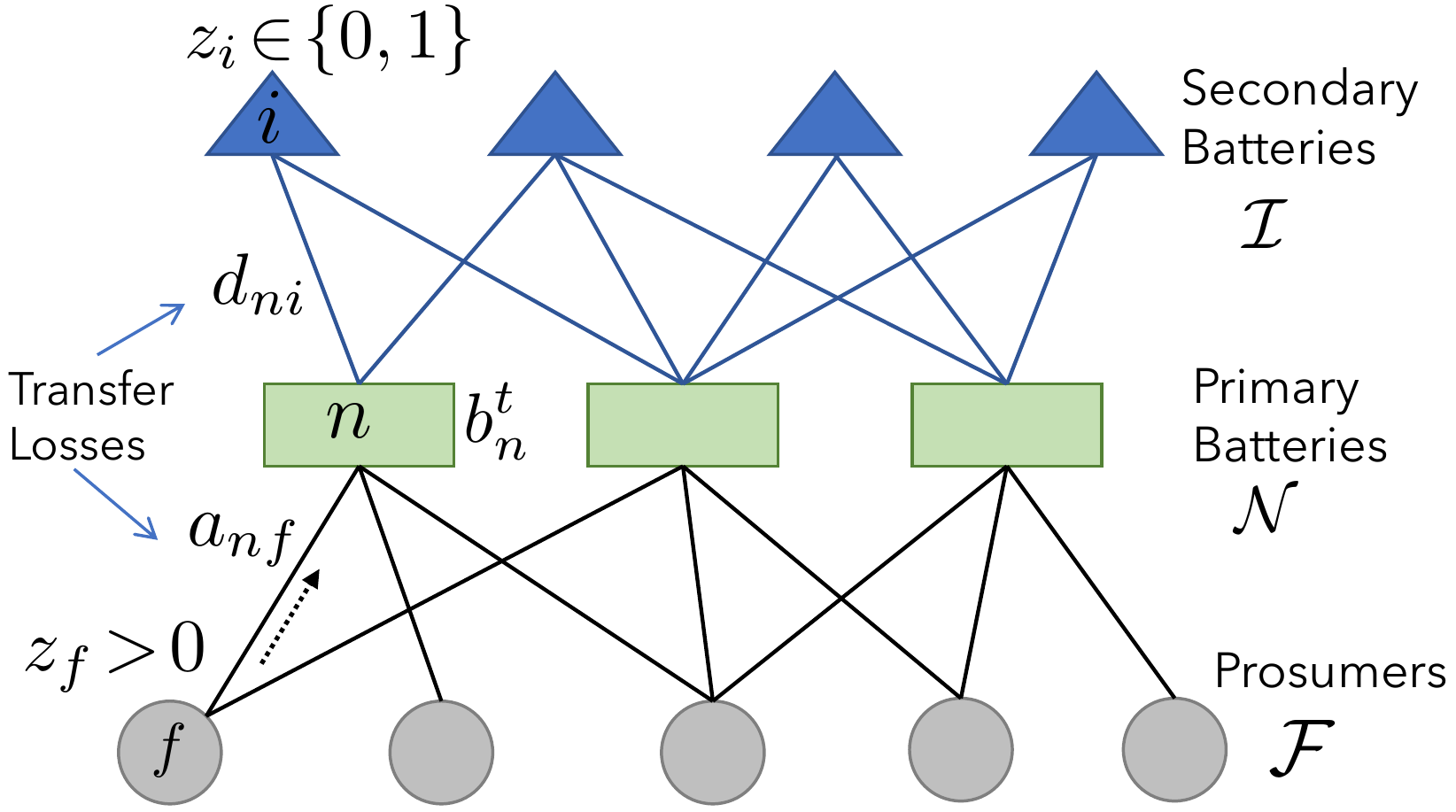}
    \caption{Microgrid energy storage optimization for a community of $\mathcal F$ prosumers and a two-tier battery system.}\label{fig:energy_problem}
\end{figure}

\textbf{Complexity}. $\mathbb{P}_{4}^t$ is NP-complete as it falls in the category of multi-dimensional Knapsack problems, see \citeAP{knapsacks-martello}.

\textbf{Datasets.} We randomly generated from a uniform distribution the scalars $a_{i,n}$ representing the resources needed to ship commodity $n$ from the $i$-th facility. Furthermore, for each time series the time-varying capacity and revenue vectors are defined as follows: 
\begin{equation*}
\begin{aligned}
{\bf{c^{t+1}}} &= {{\bf{c^t}}} + a_1 \sin(t/{\tt{period}_1})+ a_2 \sin(t/{\tt{period}_2}) + {\bf{w}}_c, \\
{\bf{b^{t+1}}} &= {{\bf{b^t}}} + {\bf{w}}_b, 
\end{aligned}
\end{equation*}
where ${\bf{c^{t}}} = [c^t_{1}, \ldots, c^t_{N}]$, ${\bf{b_{t}}} = [b^t_{1}, \ldots, b^t_{ |\mathcal{I}|}]$, the scalars $a_1$, $a_2$ are randomly generated for each time series, and the random vectors ${\bf{w}}_c \in \mathbb{R}^{|\mathcal{I}||\mathcal{N}|}$ and ${\bf{w}}_b \in \mathbb{R}^{|\mathcal{I}|}$. For further details about the number of time series and problem instances per time series please refer to Table~\ref{tab:timeseries}.

\subsubsection{Energy Grid}
We consider a microgrid energy-sharing problem, cf. \citeAP{iordanis-fnt}, where a set $\mathcal F$ of $F=|\mathcal F|$ energy producers-consumers (or, \emph{prosumers}) coordinate their energy prosumption plan over several time periods. The prosumers have at their disposal a hierarchical energy storage system, where the primary batteries are already in place and the secondary batteries are deployed upon demand by paying additional cost. The goal is to optimize the overall microgrid operation, by solving the \emph{mixed binary problem}: 
\begin{equation*}
\begin{aligned}
\mathbb{P}_{5}^t: \underset{\vz}{\text{maximize}} &\quad \sum_{f \in \mathcal{F}} c^t_{f} z^{(c)}_{f} - \sum_{i \in \mathcal{I}} p^t_{i} z^{(b)}_{i} \\
\text{subject to} & \quad \sum_{f\in\mathcal{F}} a_{nf} z^{(c)}_{f} \leq b^t_n + \sum_{i \in \mathcal{I}}d_{n, i} z^{(b)}_i, \quad \forall n \in \mathcal{N}, \\
& \quad \sum_{f \in \mathcal{F}} z^{(c)}_f = 1, \\
& \quad \vzc \in \mathbb R^{F}, \, \vzb\in \{0,1\}^F, \\
\end{aligned}
\end{equation*}
where the last constraint enforces the prosumers $\mathcal F$ to keep at least 1 (normalized) unit of energy for serving the needs of their community. 

The main \emph{variables} of $\mathbb{P}_{5}^t$ are:
\begin{itemize}
	\item $z^{(c)}_{f}\in \mathbb R$: continuous variable deciding how much energy prosumer $f\in\mathcal F$ will store (if positive) or purchase (if negative),\footnote{Negative transfers release capacity that can be used for stored-energy transfers.}
	\item $z^{(b)}_{i}\in\{0,1\}$: decides whether to deploy or not ($z_i=0$) the secondary battery $i\in\mathcal I$,
\end{itemize}
and the main \emph{parameters}:
\begin{itemize}
	\item $a_{n,f}\in \mathbb R_+$: energy transfer loss coefficient when prosumer $f\in \mathcal F$ transfers energy to battery $n\in\mathcal N$. Larger values indicate larger losses,
	\item $c_f^t\in \mathbb R_+$: selling or purchase (when $z_f^{(c)}\!<\!0$) price for the energy of prosumer $f\in \mathcal F$, at time $t$,
	\item $b_n^t\in \mathbb R_+$: available capacity at primary battery $n\in\mathcal N$ at time $t$,
	\item $d_{n,i}\in \mathbb R_+$: energy loss coefficient when transferring energy from primary battery $n\in \mathcal N$ to secondary battery $i\in\mathcal I$.
\end{itemize}

\textbf{Complexity}. $\mathbb{P}_{5}^t$ is NP-hard as it generalizes the facility location problem \citeAP{iordanis-fnt}.

 \textbf{Datasets}. We created instances using a non-stationary formula for the time-varying vectors that exhibits temporal properties as follows: 
\begin{equation*}
\begin{aligned}
\bm{c}^{t+1} &= \bm{c}^t + a \sin(t/{\tt{period}_1}) + {\bf{w}}_c, \\
\bm{b}^{t+1} &= \bm{b}^t + \bf{w}_b, 
\end{aligned}
\end{equation*}
where the scalar $a$ is randomly generated for each time series, and the random vectors ${\bf{w}}_c \in \mathbb{R}^{|\mathcal{I}||\mathcal{P}|}$ and ${\bf{w}}_b \in \mathbb{R}^{|\mathcal{I}|}$.

\subsubsection{Caching}
Our final problem is a standard data caching problem where we wish to store the most popular files in a cache of limited capacity, see \citeAP{paschos-fnt}. In particular, we consider a content library of $\mathcal N=\{1,2,\ldots,N\}$ files and a single cache with $C$ bytes capacity. Each file $n\in \mathcal N$ has popularity $p_{n}^t$ at time $t$, which captures the expected requests for file $n$ at that time, and size of $q_n$ bytes. The goal is to store, at each time, those files that will be requested by the larger number of content viewers, i.e., to maximize the cache hits. This can be formalized with the following \emph{binary problem}: 
\begin{equation*}
\begin{aligned}
\mathbb{P}_{6}^t: \underset{\vz}{\text{maximize}} &\quad \sum_{n=1}^N p_{n}^tx_{n} \\
\text{subject to} & \quad \sum_{n=1}^N x_{n} q_{n} \leq C\\
& \quad x_n\in\{0,1\}, \,\,\, \forall n\in \mathcal N.
\end{aligned}
\end{equation*}

The \emph{variables} of $\mathbb{P}_{6}^t$ are:
\begin{itemize}
	\item $x_n\in \{0,1\}$: binary variables that select a file to be cached or not,
\end{itemize}
and the main \emph{parameters}:
\begin{itemize}
	\item $p_{n}^t\in \mathbb R_+$: normalized file popularity parameters at time $t$,
	\item $q_n \in \mathbb R_+$: size of file $n\in \mathcal N$.
\end{itemize}

\textbf{Complexity}. $\mathbb{P}_{6}^t$ is a standard Knapsack problem which is NP-Complete but can be solved in pseudo-polynomial time through dynamic programming.

 \textbf{Datasets}. The file sizes were created randomly from a uniform distribution $\mathcal{U}(1,10)$, while for the file popularity we used the standard dataset for request traces from Movielens \cite{movielens}, which we sliced in order to create the different time series.

%% file: main.bbl
\begin{thebibliography}{5}
\providecommand{\natexlab}[1]{#1}
\providecommand{\url}[1]{\texttt{#1}}
\expandafter\ifx\csname urlstyle\endcsname\relax
  \providecommand{\doi}[1]{doi: #1}\else
  \providecommand{\doi}{doi: \begingroup \urlstyle{rm}\Url}\fi

\bibitem[{A. R. Karlin, N. Klein, and S. O. Gharan}(2022)]{karlin-tsp22}
{A. R. Karlin, N. Klein, and S. O. Gharan}.
\newblock {A (Slightly) Improved Approximation Algorithm for Metric TSP}.
\newblock \emph{{Arxiv}}, 2022.
\newblock URL \url{arXiv:2007.01409}.

\bibitem[{G. Dantzig, R. Fulkerson, and S. Johnson}(1954)]{dantzig-routing}
{G. Dantzig, R. Fulkerson, and S. Johnson}.
\newblock {Solution of a Large-Scale Traveling-Salesman Problem}.
\newblock \emph{{Journal of the Operations Research Society of America}},
  2\penalty0 (4):\penalty0 393--410, 1954.

\bibitem[{I. Koutsopoulos, T. G. Papaioannou, and V.
  Hatzi}(2016)]{iordanis-fnt}
{I. Koutsopoulos, T. G. Papaioannou, and V. Hatzi}.
\newblock {Modeling and Optimization of the Smart Grid Ecosystem}.
\newblock \emph{{Found. Trends Netw.}}, 10\penalty0 (2-3), 2016.

\bibitem[{Martello, Silvano and Toth, Paolo}(1990)]{knapsacks-martello}
{Martello, Silvano and Toth, Paolo}.
\newblock {Knapsack Problems: Algorithms and Computer Implementations}.
\newblock \emph{{John Wiley \& Sons}}, ISBN 978-0471924203, 1990.

\bibitem[{S. Li}(2013)]{ufl}
{S. Li}.
\newblock {A 1.488 Approximation Algorithm for the Uncapacitated Facility
  Location Problem}.
\newblock \emph{{Information and Computation}}, 222:\penalty0 45--58, 2013.

\end{thebibliography}


\begin{thebibliography}{62}
\providecommand{\natexlab}[1]{#1}
\providecommand{\url}[1]{\texttt{#1}}
\expandafter\ifx\csname urlstyle\endcsname\relax
  \providecommand{\doi}[1]{doi: #1}\else
  \providecommand{\doi}{doi: \begingroup \urlstyle{rm}\Url}\fi

\bibitem[{A. H. Land, and A. G. Doig}(1960)]{land-bnb}
{A. H. Land, and A. G. Doig}.
\newblock {An Automatic Method of Solving Discrete Programming Problems}.
\newblock \emph{Econometrica}, 28\penalty0 (3):\penalty0 497--520, 1960.

\bibitem[{A. M. Alvarez}(2017)]{alvarez-branching-17}
{A. M. Alvarez}.
\newblock {A Machine Learning-based Approximation of Strong Branching}.
\newblock \emph{{INFORMS Journal on Computing}}, 29\penalty0 (1):\penalty0
  185--195, 2017.

\bibitem[Agrawal et~al.(2019)Agrawal, Amos, Barratt, Boyd, Diamond, and
  Kolter]{agrawal2019differentiable}
Akshay Agrawal, Brandon Amos, Shane Barratt, Stephen Boyd, Steven Diamond, and
  J~Zico Kolter.
\newblock Differentiable convex optimization layers.
\newblock \emph{Advances in Neural Information Processing Systems (NeurIPS)},
  2019.

\bibitem[Ans{\'o}tegui et~al.(2015)Ans{\'o}tegui, Malitsky, Samulowitz,
  Sellmann, and Tierney]{ansotegui-aistats15}
Carlos Ans{\'o}tegui, Yuri Malitsky, Horst Samulowitz, Meinolf Sellmann, and
  Kevin Tierney.
\newblock Model-based genetic algorithms for algorithm configuration.
\newblock In \emph{Twenty-Fourth International Joint Conference on Artificial
  Intelligence}, 2015.

\bibitem[Ba et~al.(2016)Ba, Kiros, and Hinton]{ba2016layer}
Jimmy~Lei Ba, Jamie~Ryan Kiros, and Geoffrey~E Hinton.
\newblock Layer normalization.
\newblock \emph{arXiv preprint arXiv:1607.06450}, 2016.

\bibitem[Balcan et~al.(2018)Balcan, Dick, Sandholm, and
  Vitercik]{balcan-branch-18}
Maria-Florina Balcan, Travis Dick, Tuomas Sandholm, and Ellen Vitercik.
\newblock Learning to branch.
\newblock In \emph{International Conference on Machine Learning (ICML)}, pages
  344--353, 2018.

\bibitem[Balcan et~al.(2021)Balcan, Prasad, Sandholm, and
  Vitercik]{balcan2021sample}
Maria-Florina~F Balcan, Siddharth Prasad, Tuomas Sandholm, and Ellen Vitercik.
\newblock Sample complexity of tree search configuration: Cutting planes and
  beyond.
\newblock \emph{Advances in Neural Information Processing Systems (NeurIPS)},
  2021.

\bibitem[Battaglia et~al.(2018)Battaglia, Hamrick, Bapst, Sanchez-Gonzalez,
  Zambaldi, Malinowski, Tacchetti, Raposo, Santoro, Faulkner,
  et~al.]{battaglia2018relational}
Peter~W Battaglia, Jessica~B Hamrick, Victor Bapst, Alvaro Sanchez-Gonzalez,
  Vinicius Zambaldi, Mateusz Malinowski, Andrea Tacchetti, David Raposo, Adam
  Santoro, Ryan Faulkner, et~al.
\newblock Relational inductive biases, deep learning, and graph networks.
\newblock \emph{arXiv preprint arXiv:1806.01261}, 2018.

\bibitem[Bengio et~al.(2021)Bengio, Lodi, and Prouvost]{bengio2021machine}
Yoshua Bengio, Andrea Lodi, and Antoine Prouvost.
\newblock Machine learning for combinatorial optimization: a methodological
  tour d’horizon.
\newblock \emph{European Journal of Operational Research}, 290\penalty0
  (2):\penalty0 405--421, 2021.

\bibitem[Bertsimas and Stellato(2019)]{bertsimas2019online}
Dimitris Bertsimas and Bartolomeo Stellato.
\newblock {Online Mixed-integer Optimization in Milliseconds}.
\newblock \emph{arXiv preprint arXiv:1907.02206}, 2019.

\bibitem[Bertsimas and Stellato(2021)]{bertsimas2021voice}
Dimitris Bertsimas and Bartolomeo Stellato.
\newblock {The Voice of Optimization}.
\newblock \emph{Machine Learning}, 110\penalty0 (2):\penalty0 249--277, 2021.

\bibitem[Bonami et~al.(2018)Bonami, Lodi, and Zarpellon]{bonami18}
Pierre Bonami, Andrea Lodi, and Giulia Zarpellon.
\newblock Learning a classification of mixed-integer quadratic programming
  problems.
\newblock In \emph{International Conference on the Integration of Constraint
  Programming, Artificial Intelligence, and Operations Research}, pages
  595--604. Springer, 2018.

\bibitem[Borrelli et~al.(2017)Borrelli, Bemporad, and
  Morari]{borrelli2017predictive}
Francesco Borrelli, Alberto Bemporad, and Manfred Morari.
\newblock \emph{Predictive control for linear and hybrid systems}.
\newblock Cambridge University Press, 2017.

\bibitem[Chen(2010)]{chen10}
Zhi-Long Chen.
\newblock Integrated production and outbound distribution scheduling: review
  and extensions.
\newblock \emph{Operations research}, 58\penalty0 (1):\penalty0 130--148, 2010.

\bibitem[Clenshaw and Curtis(1960)]{clenshaw1960method}
Charles~W Clenshaw and Alan~R Curtis.
\newblock A method for numerical integration on an automatic computer.
\newblock \emph{Numerische Mathematik}, 2\penalty0 (1):\penalty0 197--205,
  1960.

\bibitem[Dai et~al.(2016)Dai, Dai, and Song]{dai-icml16}
Hanjun Dai, Bo~Dai, and Le~Song.
\newblock Discriminative embeddings of latent variable models for structured
  data.
\newblock In \emph{International conference on machine learning}, pages
  2702--2711, 2016.

\bibitem[Dilkina et~al.(2009{\natexlab{a}})Dilkina, Gomes, Malitsky, Sabharwal,
  and Sellmann]{dilkina09-16}
Bistra Dilkina, Carla~P Gomes, Yuri Malitsky, Ashish Sabharwal, and Meinolf
  Sellmann.
\newblock Backdoors to combinatorial optimization: Feasibility and optimality.
\newblock In \emph{International Conference on Integration of Constraint
  Programming, Artificial Intelligence, and Operations Research}, pages 56--70.
  Springer, 2009{\natexlab{a}}.

\bibitem[Dilkina et~al.(2009{\natexlab{b}})Dilkina, Gomes, and
  Sabharwal]{dilkina09-17}
Bistra Dilkina, Carla~P Gomes, and Ashish Sabharwal.
\newblock Backdoors in the context of learning.
\newblock In \emph{International Conference on Theory and Applications of
  Satisfiability Testing}, pages 73--79. Springer, 2009{\natexlab{b}}.

\bibitem[Ding et~al.(2020)Ding, Zhang, Shen, Li, Wang, Xu, and
  Song]{dingaaai20}
Jian-Ya Ding, Chao Zhang, Lei Shen, Shengyin Li, Bing Wang, Yinghui Xu, and
  Le~Song.
\newblock Accelerating primal solution findings for mixed integer programs
  based on solution prediction.
\newblock In \emph{Proceedings of the AAAI Conference on Artificial
  Intelligence}, volume~34, pages 1452--1459, 2020.

\bibitem[{E. B. Khalil, B. Dilkina, G. L. Nemhauser, S. Ahmed, and Y.
  Shao}(2017)]{khalil-branch-17}
{E. B. Khalil, B. Dilkina, G. L. Nemhauser, S. Ahmed, and Y. Shao}.
\newblock {Learning to Run Heuristics in Tree Search}.
\newblock In \emph{International Joint Conference on Artificial Intelligence
  (IJCAI)}, 2017.

\bibitem[{E. B. Khalil, et al.}(2016)]{khalil-branch-16}
{E. B. Khalil, et al.}
\newblock {Learning to Branch in Mixed Integer Programming}.
\newblock In \emph{Proceedings of the AAAI conference on artificial
  intelligence}, page 724–731, 2016.

\bibitem[{F. Hutter, et al.}(2009)]{hutter09}
{F. Hutter, et al.}
\newblock {ParamILS: An Automatic Algorithm Configuration Framework}.
\newblock \emph{{Journal of Artificial Intelligence Research}}, 36\penalty0
  (1), 2009.

\bibitem[{F. Hutter, et al.}(2011)]{hutter11}
{F. Hutter, et al.}
\newblock {Sequential Model-based Optimization for General Algorithm
  Configuration}.
\newblock \emph{{Learning and intelligent optimization}}, pages 507--523, 2011.

\bibitem[{F. M. Harper, and J. A. Konstan}(2015)]{movielens}
{F. M. Harper, and J. A. Konstan}.
\newblock {The MovieLens Datasets: History and Context}.
\newblock \emph{{ACM Transactions on Interactive Intelligent Systems}},
  5:\penalty0 19:1--19:2, 2015.

\bibitem[Farahani and Hekmatfar(2009)]{farahani09}
R.~Z. Farahani and M.~Hekmatfar.
\newblock {Facility Location: Concepts, Models, Algorithms and Case Studies}.
\newblock \emph{{Springer}}, 2009.

\bibitem[Faísca et~al.(2007)Faísca, Dua, and Pistikopoulos]{MPLP}
Nuno~P. Faísca, Vivek Dua, and Efstratios~N. Pistikopoulos.
\newblock \emph{Multiparametric Linear and Quadratic Programming}, chapter~1,
  pages 1--23.
\newblock John Wiley \& Sons, Ltd, 2007.
\newblock ISBN 9783527631216.
\newblock \doi{https://doi.org/10.1002/9783527631216.ch1}.
\newblock URL
  \url{https://onlinelibrary.wiley.com/doi/abs/10.1002/9783527631216.ch1}.

\bibitem[{G. Barlacchi, et al.}(2015)]{torino-dataset}
{G. Barlacchi, et al.}
\newblock {A Multi-source Dataset of Urban Life in the City of Milan and the
  Province of Trentino}.
\newblock \emph{{Scientific Data}}, 2, 2015.

\bibitem[{G. Laporte}(2009)]{laporte09}
{G. Laporte}.
\newblock {Fifty Years of Vehicle Routing}.
\newblock \emph{{Transportation Science}}, pages 408--416, 2009.

\bibitem[{G. S. Paschos, et al.}(2019)]{paschos-fnt}
{G. S. Paschos, et al.}
\newblock {Cache Optimization Models and Algorithms}.
\newblock \emph{{Found. Trends Commun. Inf. Theory}}, 16\penalty0
  (3-4):\penalty0 156--343, 2019.

\bibitem[{G. Vivek, K. Srikumar, and D. Roth}(2012)]{vivek2012}
{G. Vivek, K. Srikumar, and D. Roth}.
\newblock {On amortizing inference cost for structured prediction}.
\newblock \emph{{Joint Conference on Empirical Methods in Natural Language
  Processing and Computational Natural Language Learning}}, pages 1114--1124,
  2012.

\bibitem[Gasse et~al.(2019)Gasse, Ch{\'e}telat, Ferroni, Charlin, and
  Lodi]{gasse2019exact}
Maxime Gasse, Didier Ch{\'e}telat, Nicola Ferroni, Laurent Charlin, and Andrea
  Lodi.
\newblock Exact combinatorial optimization with graph convolutional neural
  networks.
\newblock \emph{Advances in Neural Information Processing Systems (NeurIPS)},
  2019.

\bibitem[Geifman and El-Yaniv(2019)]{geifman2019selectivenet}
Yonatan Geifman and Ran El-Yaniv.
\newblock Selectivenet: A deep neural network with an integrated reject option.
\newblock In \emph{International Conference on Machine Learning (ICML)}, pages
  2151--2159, 2019.

\bibitem[Gleixner et~al.(2021)Gleixner, Hendel, Gamrath, Achterberg, Bastubbe,
  Berthold, Christophel, Jarck, Koch, Linderoth, L\"ubbecke, Mittelmann,
  Ozyurt, Ralphs, Salvagnin, and Shinano]{miplib}
Ambros Gleixner, Gregor Hendel, Gerald Gamrath, Tobias Achterberg, Michael
  Bastubbe, Timo Berthold, Philipp~M. Christophel, Kati Jarck, Thorsten Koch,
  Jeff Linderoth, Marco L\"ubbecke, Hans~D. Mittelmann, Derya Ozyurt, Ted~K.
  Ralphs, Domenico Salvagnin, and Yuji Shinano.
\newblock {MIPLIB 2017: Data-Driven Compilation of the 6th Mixed-Integer
  Programming Library}.
\newblock \emph{Mathematical Programming Computation}, 2021.
\newblock \doi{10.1007/s12532-020-00194-3}.
\newblock URL \url{https://doi.org/10.1007/s12532-020-00194-3}.

\bibitem[Gopalakrishnan and Johnson(2005)]{johnson05}
B.~Gopalakrishnan and E.~L. Johnson.
\newblock {Airline Crew Scheduling: State-of-the-Art}.
\newblock \emph{{Annals of Operations Research}}, 140\penalty0 (1):\penalty0
  305--337, 2005.

\bibitem[Hansknecht et~al.(2018)Hansknecht, Joormann, and
  Stiller]{hansknecht-branch18}
Christoph Hansknecht, Imke Joormann, and Sebastian Stiller.
\newblock Cuts, primal heuristics, and learning to branch for the
  time-dependent traveling salesman problem.
\newblock \emph{arXiv preprint arXiv:1805.01415}, 2018.

\bibitem[He et~al.(2014)He, Daume~III, and Eisner]{he-neurips14}
He~He, Hal Daume~III, and Jason~M Eisner.
\newblock Learning to search in branch and bound algorithms.
\newblock \emph{Advances in Neural Information Processing Systems (NeurIPS)},
  2014.

\bibitem[{I. Boussaıd, et al.}(2013)]{boussaid-heuristics}
{I. Boussaıd, et al.}
\newblock {A Survey on Optimization Metaheuristics}.
\newblock \emph{{Information Sciences}}, 237:\penalty0 82--117, 2013.

\bibitem[{J. Song, R. Lanka, A. Zhao, Y. Yue, and M. Ono}(2018)]{song-arxiv18}
{J. Song, R. Lanka, A. Zhao, Y. Yue, and M. Ono}.
\newblock {Learning to Search via Retrospective Imitation}.
\newblock In \emph{arXiv:1804.00846}, 2018.

\bibitem[Karalias and Loukas(2020)]{karalias202Nips}
Nikolaos Karalias and Andres Loukas.
\newblock Erdos goes neural: an unsupervised learning framework for
  combinatorial optimization on graphs.
\newblock In \emph{Proceedings of NeurIPS}, 2020.

\bibitem[Khalil et~al.(2017)Khalil, Dai, Zhang, Dilkina, and
  Song]{dai2017learning}
Elias Khalil, Hanjun Dai, Yuyu Zhang, Bistra Dilkina, and Le~Song.
\newblock Learning combinatorial optimization algorithms over graphs.
\newblock \emph{Advances in Neural Information Processing Systems (NeurIPS)},
  2017.

\bibitem[Kim et~al.(2021)Kim, Park, et~al.]{kim2021learning}
Minsu Kim, Jinkyoo Park, et~al.
\newblock Learning collaborative policies to solve np-hard routing problems.
\newblock \emph{Advances in Neural Information Processing Systems (NeurIPS)},
  2021.

\bibitem[Kipf and Welling(2016)]{kipf2016semi}
Thomas~N Kipf and Max Welling.
\newblock Semi-supervised classification with graph convolutional networks.
\newblock \emph{arXiv preprint arXiv:1609.02907}, 2016.

\bibitem[Klau{\v{c}}o et~al.(2019)Klau{\v{c}}o, Kal{\'u}z, and
  Kvasnica]{klauco2019}
Martin Klau{\v{c}}o, Martin Kal{\'u}z, and Michal Kvasnica.
\newblock Machine learning-based warm starting of active set methods in
  embedded model predictive control.
\newblock \emph{Engineering Applications of Artificial Intelligence},
  77:\penalty0 1--8, 2019.

\bibitem[Kleinberg et~al.(2019)Kleinberg, Leyton-Brown, Lucier, and
  Graham]{kleinberg-neurips19}
Robert Kleinberg, Kevin Leyton-Brown, Brendan Lucier, and Devon Graham.
\newblock Procrastinating with confidence: Near-optimal, anytime, adaptive
  algorithm configuration.
\newblock \emph{Advances in Neural Information Processing Systems (NeurIPS)},
  2019.

\bibitem[Kruber et~al.(2017)Kruber, L{\"u}bbecke, and Parmentier]{kruber17}
Markus Kruber, Marco~E L{\"u}bbecke, and Axel Parmentier.
\newblock Learning when to use a decomposition.
\newblock In \emph{International conference on AI and OR techniques in
  constraint programming for combinatorial optimization problems}, pages
  202--210. Springer, 2017.

\bibitem[Kulesza and Pereira(2007)]{kulesza08}
Alex Kulesza and Fernando Pereira.
\newblock Structured learning with approximate inference.
\newblock \emph{Advances in Neural Information Processing Systems (NeurIPS)},
  2007.

\bibitem[Lakhina et~al.(2004)Lakhina, Papagiannaki, Crovella, Diot, Kolaczyk,
  and Taft]{diotsigm04}
Anukool Lakhina, Konstantina Papagiannaki, Mark Crovella, Christophe Diot,
  Eric~D. Kolaczyk, and Nina Taft.
\newblock Structural analysis of network traffic flows.
\newblock In \emph{Proceedings of ACM Sigmetrics}, 2004.

\bibitem[Leng et~al.(2018)Leng, Dou, Li, Zhu, and Jin]{leng-aaai18}
Cong Leng, Zesheng Dou, Hao Li, Shenghuo Zhu, and Rong Jin.
\newblock Extremely low bit neural network: Squeeze the last bit out with
  {ADMM}.
\newblock In \emph{Proceedings of the AAAI Conference on Artificial
  Intelligence}, 2018.

\bibitem[{M. C. Gonzalez, et al.}(2008)]{hidalgo-2008}
{M. C. Gonzalez, et al.}
\newblock {Understanding Individual Human Mobility Patterns}.
\newblock \emph{{Nature}}, 453:\penalty0 779--782, 2008.

\bibitem[{M. Lopez-Ibanez, et al.}(2016)]{lopez16}
{M. Lopez-Ibanez, et al.}
\newblock {The Irace Package: Iterated Racing for Automatic Algorithm
  Configuration}.
\newblock \emph{{Operations Research Perspectives}}, pages 43--58, 2016.

\bibitem[Misra et~al.(2022)Misra, Roald, and Ng]{misra2019}
Sidhant Misra, Line Roald, and Yeesian Ng.
\newblock Learning for constrained optimization: Identifying optimal active
  constraint sets.
\newblock \emph{INFORMS Journal on Computing}, 34\penalty0 (1):\penalty0
  463--480, 2022.

\bibitem[Nair et~al.(2020)Nair, Bartunov, Gimeno, von Glehn, Lichocki, Lobov,
  O'Donoghue, Sonnerat, Tjandraatmadja, Wang, et~al.]{nair2020solving}
Vinod Nair, Sergey Bartunov, Felix Gimeno, Ingrid von Glehn, Pawel Lichocki,
  Ivan Lobov, Brendan O'Donoghue, Nicolas Sonnerat, Christian Tjandraatmadja,
  Pengming Wang, et~al.
\newblock Solving mixed integer programs using neural networks.
\newblock \emph{arXiv preprint arXiv:2012.13349}, 2020.

\bibitem[Orlowski et~al.()Orlowski, Pi{\'o}ro, Tomaszewski, and
  Wess{\"a}ly]{SNDlib10}
S.~Orlowski, M.~Pi{\'o}ro, A.~Tomaszewski, and R.~Wess{\"a}ly.
\newblock {SNDlib} 1.0--{S}urvivable {N}etwork {D}esign {L}ibrary.
\newblock In \emph{Proceedings of the 3rd International Network Optimization
  Conference (INOC 2007), Spa, Belgium}, April .
\newblock http://sndlib.zib.de, extended version accepted in Networks, 2009.

\bibitem[{R. E. Bixby}(2010)]{bixby}
{R. E. Bixby}.
\newblock {A Brief History of Linear and Mixed-Integer Programming
  Compucation}.
\newblock \emph{{Documenta Mathematica}}, pages 107--121, 2010.

\bibitem[{S. Chopra, I. Gilboa, and S. T. Sastry}(1998)]{routing-complexity}
{S. Chopra, I. Gilboa, and S. T. Sastry}.
\newblock {Source Sink Flows with Capacity Installation in Batches}.
\newblock \emph{{Discrete Applied Mathematics}}, 85:\penalty0 165--192, 1998.

\bibitem[{S. Hochreiter, et al.}(1997)]{lstm97}
{S. Hochreiter, et al.}
\newblock Long short-term memory.
\newblock \emph{{Neural Computation}}, 9\penalty0 (8):\penalty0 1735--1780,
  1997.

\bibitem[Song et~al.(2020)Song, Yue, Dilkina, et~al.]{song2020general}
Jialin Song, Yisong Yue, Bistra Dilkina, et~al.
\newblock A general large neighborhood search framework for solving integer
  linear programs.
\newblock \emph{Advances in Neural Information Processing Systems (NeurIPS)},
  2020.

\bibitem[{T. Achterberg}(2009)]{scip}
{T. Achterberg}.
\newblock {SCIP: Solving Constraint Integer Programs}.
\newblock \emph{{Mathematical Programming Computation}}, 1:\penalty0 1--41,
  2009.

\bibitem[Williamson and Shmoys(2011)]{shmoys-book}
David Williamson and David Shmoys.
\newblock \emph{{The Design of Approximation Algorithms}}.
\newblock Cambridge University Press, 2011.

\bibitem[Wu et~al.(2021)Wu, Song, Cao, and Zhang]{wu2021learning}
Yaoxin Wu, Wen Song, Zhiguang Cao, and Jie Zhang.
\newblock Learning large neighborhood search policy for integer programming.
\newblock \emph{Advances in Neural Information Processing Systems (NeurIPS)},
  2021.

\bibitem[Xavier et~al.(2021)Xavier, Qiu, and Ahmed]{xavier-informs21}
Alinson~S Xavier, Feng Qiu, and Shabbir Ahmed.
\newblock Learning to solve large-scale security-constrained unit commitment
  problems.
\newblock \emph{INFORMS Journal on Computing}, 33\penalty0 (2):\penalty0
  739--756, 2021.

\bibitem[{Y. Pochet, and L. A. Wolsey}(2006)]{pochet2006}
{Y. Pochet, and L. A. Wolsey}.
\newblock {Production Planning by Mixed Integer Programming}.
\newblock \emph{{ Springer Science and Business Media}}, 2006.

\end{thebibliography}
